\documentclass[opre,nonblindrev]{informs3} % current default for manuscript submission

\DoubleSpacedXI % Made default 4/4/2014 at request
\usepackage{endnotes}
\let\footnote=\endnote
\let\enotesize=\normalsize
\def\notesname{Endnotes}%
\def\makeenmark{$^{\theenmark}$}
\def\enoteformat{\rightskip0pt\leftskip0pt\parindent=1.75em
  \leavevmode\llap{\theenmark.\enskip}}

\usepackage{natbib}
 \bibpunct[, ]{(}{)}{,}{a}{}{,}%
 \def\bibfont{\small}%

\usepackage{lscape} 
\usepackage{rotating}
\usepackage{graphicx}

\TheoremsNumberedThrough     % Preferred (Theorem 1, Lemma 1, Theorem 2)
\ECRepeatTheorems

\EquationsNumberedThrough    % Default: (1), (2), ...
\begin{document}
%%%%%%%%%%%%%%%%

%\setcounter{page}{1}
%\DOI{10.1287/xxxx.0000.0000}%

\RUNAUTHOR{Ton de Kok and Mirjam Meijer}

\RUNTITLE{Fixed non-stockout-probability policies for the single-item lost-sales model}

% \TITLE{Bundling Information Goods of Decreasing Value}
\TITLE{Fixed non-stockout-probability policies for the single-item lost-sales model}

\ARTICLEAUTHORS{%
\AUTHOR{Ton de Kok}
\AFF{School of Industrial Engineering, Eindhoven University of Technology, P.O. Box 513, 5600MB Eindhoven, The Netherlands \EMAIL{a.g.d.kok@tue.nl}} %, \URL{}}
\AUTHOR{Mirjam S. Meijer}
\AFF{Kühne Logistics University,  Großer Grasbrook 17, 20457 Hamburg, Germany 
\EMAIL{mirjam.meijer@klu.org}} %, \URL{}}
} % end of the block

\ABSTRACT{%
\textbf{Abstract}
\newline
We consider the classical discrete time lost-sales model under stationary continuous demand and linear holding and penalty costs and positive constant lead time. To date the optimal policy structure is only known implicitly by solving numerically the Bellman equations. In this paper we derive an optimality equation for the lost-sales model. We propose a fixed non-stockout-probability ($FP_3$) policy, implying that each period the order size ensures that $P_3$, the probability of no-stockout at the end of the period of arrival of this order, equals some target value. The $FP_3$-policy can be computed efficiently and accurately from an exact recursive expression and two-moment fits to the emerging random variables. We use the lost-sales optimality equation to compute the optimal $FP_3$-policy. Comparison against the optimal policy for discrete demand suggests that the fixed $P_3$-policy is close-to-optimal. An extensive numerical experiment shows that the $FP_3$-policy outperforms other policies proposed in literature in $97\%$ of all cases. Under the $FP_3$-policy, the volatility of the replenishment process, measured as coefficient of variation (cv) is much lower than the volatility of the demand process. This cv-reduction holds a promise for substantial cost reduction at upstream stages in the supply chain of the end-item under consideration, compared to the situation with backlogging.\\

\textbf{Subject classifications}
\newline
inventory/production: stochastic lost-sales model, review/lead times, optimality equation; 
probability: conditional random variables
simulation: simulation-based optimization
% Enter your abstract
}%

%\KEYWORDS{deterministic inventory theory; infinite linear programming duality;
%  existence of optimal policies; semi-Markov decision process; cyclic schedule}

\KEYWORDS{lost-sales model, fixed non-stockout probability policy, conditional random variables, moment-iteration} 
\HISTORY{This paper was
first submitted on April 8, 2023. A revision was submitted on October 15, 2024}

\maketitle

%%%%%%%%%%%%%%%%%%%%%%%%%%%%%%%%%%%%%%%%%%%%%%%%%%%%%%%%%%%%%%%%%%%%%%

% Samples of sectioning (and labeling) in OPRE
% NOTE: (1) \section and \subsection do NOT end with a period
%       (2) \subsubsection and lower need end punctuation
%       (3) capitalization is as shown (title style).
%
%\section{Introduction.}\label{intro} %%1.
%\subsection{Duality and the Classical EOQ Problem.}\label{class-EOQ} %% 1.1.
%\subsection{Outline.}\label{outline1} %% 1.2.
%\subsubsection{Cyclic Schedules for the General Deterministic SMDP.}
%  \label{cyclic-schedules} %% 1.2.1
%\section{Problem Description.}\label{problemdescription} %% 2.

% Text of your paper here

\section{Introduction}
\label{sec:intro}
We consider the classical periodic review single-item-single-echelon lost-sales model with constant integer lead times and linear holding and penalty costs. Over the last decade this model received considerable interest as new results have been derived regarding the optimal policy structure. We refer to \cite{goldberg2016asymptotic}, \cite{huh2009asymptotic} for the asymptotic optimality of constant order ($CO$) policies for lead time to infinity, and the optimality of base stock ($BS$) policies for penalty costs to infinity, respectively. In this paper we investigate a particular class of policies, which was proposed by \cite{Dons1996}, who show that a myopic fixed non-stockout probability policy outperforms the base stock policy. This policy is optimal for the backlogging model. We denote this policy as $FP_3$ policy, the Fixed $P_3$ policy, using the notation in \cite{silver2016inventory}, who denote the non-stockout probability as $P_3$. Through an extensive simulation study we show that this policy outperforms the BS policy, the CO policy, the capped base stock ($CBS$) policy proposed in \cite{xin2021understanding}, and the Projected Inventory Level ($PIL$) policy proposed in \cite{van2024projected}. We also show that the costs of the $FP_3$-policy are remarkably close to the optimal costs for all cases in \cite{zipkin2008old}.

Finding the optimal $FP_3$-policy is feasible by using an optimality equation for the lost-sales model, which is presented in this paper. We argue that we can always find a target $P_3$ value for which this optimality equation is satisfied. We use discrete event simulation combined with bisection to solve for the optimal $P_3$. We derive exact expressions for $P_3$ given the state of the system, i.e. the available stock and all individual outstanding orders, and a recursive exact expression allows for an efficient approximation procedure based on two-moment fits. We apply the exact procedure proposed by \cite{van2024projected} to compute replenishment quantities under the $PIL$-policy to compute replenishment quantities under the $FP_3$-policy. We extend this exact procedure to the case of shifted mixtures of Erlang distributions with a single rate parameter. Assuming shifted exponential demand, we derive closed-form expressions for the optimal order quantity, the long-run costs, and the resulting $P_3$ under a constant order policy. These closed-form expressions are used to calibrate the simulation-based optimization methodology.  

An important auxiliary outcome of our experiments is that the $FP_3$-policy is not only cost-optimal, but also gives the lowest replenishment order coefficient of variation among all well-performing policies. This reverse bullwhip effect of the $FP_3$-policy should yield substantial cost reductions in the supply chain upstream of the stockpoint under consideration, compared with the situation where excess demand is backlogged.

The paper is structured as follows. We introduce the model and its notation in section \ref{sec2}. The lost-sales optimality equation is presented in section \ref{sec3}. In section \ref{sec4} we develop exact expressions for the targets of two myopic policies, the $FP_3$-policy and the $PIL$-policy, that must be met by the replenishment decision each period. In section \ref{Compcons} we use the exact expressions to derive efficient exact and (in most cases) accurate approximative procedures from which the replenishment quantity can be computed each period, taking into account the current inventory and all outstanding orders. An extensive numerical experiment in section \ref{sec5} underpins the above statements concerning the performance of the $FP_3$-policy compared to other policies. This motivates the generation of the optimal $P_3$ target values for a wide range of model parameters that can be used in practice together with an accurate and efficient procedure to compute the replenishment quantity each period under the $FP_3$-policy. We summarize our findings in section \ref{sec6}, where we also briefly discuss further research. 

\section{Model}
\label{sec2}
We consider the periodic review lost-sales model, where at the start of each time unit an order is released, possibly of zero size. The lead time of an order is assumed to be constant and an integer number of time units. Inventory is replenished at the start of each time unit, before a new order is released. Demand in different time units is considered i.i.d. and assumed to be continuous. Demand not satisfied immediately from stock is lost. At the end of each time unit holding costs and penalty costs are incurred. Holding costs are incurred for each unit on stock at the end of the time unit. Penalty costs for each unit lost during the time unit are incurred at the end of the time unit. In table \ref{tab: notation} we provide the notation used throughout the paper. In case we consider a steady-state equivalent of one of the random variables defined, we drop the index $t$. We use the convention that period $t$ starts at time $t-1$ and ends at time $t$. The objective of the paper is to analyze the cost performance of various policies proposed in literature and a policy we propose in this paper. We aim to minimize the long-run average cost, i.e. we do not consider a discount factor, which is given by
\begin{equation}
C = hE[\tilde{I}]+p(E[D]-E[Q]).
\end{equation}

It is well-known that to-date the optimal policy for this lost-sales model is unknown and likely to be of a complicated structure. Below we show that a property of the optimal policy can be related to the first two moments of the stationary time between two stockouts. We argue that the same property holds for two heuristic policies. Stockout moments decouple to some extent the future evolution of the system from its past, but stockout moments are not regeneration points, as the state of the system, expressed by all individual outstanding orders, is a consequence of the past. Still, the evolution of the inventory during the intervals between stockout moments allows us to derive exact expressions for the cost difference between the optimal policy and a policy that follows from perturbation of the optimal policy.\\

%In table 1 we present the notation used throughout the paper.\\

\begin{table}[ht!]
    \centering
    \begin{tabular}{ll}\\
    
    \hline
        $p$ & penalty cost per unit lost\\
        $h$ & holding cost per unit on stock at the end of a period\\
        $L$ & replenishment lead time\\
        $D_t$ & demand in period $t$\\
        $Q_t$ & quantity ordered at the start of period $t$\\
        $I_t$ & physical inventory at the start of period $t$\\
        $\tilde{I}_t$ & physical inventory at the end of period $t$\\
        $P_{3,t+L}$ & non-stockout probability at the end of period $t+L$,\\ 
        $ $ & determined at the start of period t\\
        $T_n$ & time between $(n-1)^{th}$ and $n^{th}$ stockout after time 0\\
        $T$ & stationary time between stockouts\\
        \hline
    \end{tabular}
    \caption{Notation}
    \label{tab: notation}
\end{table}

Below we omit the subscript $t$ when we consider the stationary version of the random variables defined above.\\

\section{Optimal policy}
\label{sec3}
The optimal policy for the lost-sales model is to-date unknown. In this section we derive a property of the optimal policy that is of use when determining optimal policies within classes of policies for the lost-sales model. For this purpose, we will first introduce a perturbed policy and prove some properties of this policy.

We assume that the optimal policy yields optimal actions $\{Q^{*}_t\}_{t=1}^{\infty}$. We define an $\epsilon$-perturbed optimal policy by its actions,
\begin{equation*}
\begin{aligned}
{\tilde{Q}_t(\epsilon)} &=Q^{*}_t+\epsilon, \epsilon>0.
\end{aligned}
\end{equation*}
%\\
%Note that this perturbation may imply that the optimal policy and the perturbed policy do not belong to the same class of policies. \textbf{\textcolor{blue}{Is dit een probleem?}}\textbf{\textcolor{red}{Nee, misschien weglaten. Komt doordat ik me realiseerde dat een base stock policy zo niet kan worden aagepast. Begin je met een base stock policy, dan impliceert verhogen met eps, dat de base stock policy deze eps dan de volgende periode weer "verwijderd". Het blijft natuurlijk wel een policy}}
Without loss of generality we assume that at time 0 there is a stockout under both policies. 
We define the $n^{th}$ stockout moment under the $\epsilon$-perturbed policy $S_n(\epsilon), n\geq0$, as
\begin{equation*}
\begin{aligned}
& S_0 = 0\\
& S_n (\epsilon)= \sum_{k=1}^nT_k(\epsilon), n\geq1.\\
\end{aligned}
\end{equation*}
where $T_k(\epsilon)$ denotes the time between the $(k-1)^{th}$ and $k^{th}$ stockout of the $\epsilon$-perturbed policy. Lemma \ref{lem:time_between_stockouts_increasing_in_eps} states that the time between two stockouts is increasing in $\epsilon$.\\

\begin{lemma}\label{lem:time_between_stockouts_increasing_in_eps}
    The time between two stockouts under the $\epsilon$-perturbed policy is increasing in $\epsilon$.
\end{lemma}

% \begin{proof}
\proof{Proof.} 
Take two values $\epsilon_1$ and $\epsilon_2$, such that $\epsilon_1<\epsilon_2$. Let $t=\bar{t}$ be a stockout moment and let $T(\epsilon_1)$ and $T(\epsilon_2)$ be the time until the next stockout moment under the $\epsilon_1$- and $\epsilon_2$-perturbed policies, respectively.
Since $\bar{t}+T(\epsilon_1)$ is the next stockout moment under the $\epsilon_1$-perturbed policy, it must hold that:
\begin{align*}
    D_{\bar{t}+1} &< Q^*_{\bar{t}+1-L}+\epsilon_1\\
    D_{\bar{t}+2} &< Q^*_{\bar{t}+2-L}+\epsilon_1 + Q^*_{\bar{t}+1-L}+\epsilon_1 - D_{\bar{t}+1}\\
    D_{\bar{t}+3} &< Q^*_{\bar{t}+3-L}+\epsilon_1 + Q^*_{\bar{t}+2-L}+\epsilon_1 + Q^*_{\bar{t}+1-L}+\epsilon_1 - D_{\bar{t}+1} - D_{\bar{t}+2}\\
%    &= Q_{\bar{t}+3-L}+\epsilon_2 +\sum_{i=1}^2  Q_{\bar{t}+i-L}+\epsilon_2 - D_{\bar{t}+i}\\
    \vdots & \\
    D_{\bar{t}+T(\epsilon_1)-1} &< Q^*_{\bar{t}+T(\epsilon_1)-1-L}+\epsilon_1 +\sum_{i=1}^{T(\epsilon_1)-2}  \left( Q^*_{\bar{t}+i-L}+\epsilon_1 - D_{\bar{t}+i} \right)\\
    D_{\bar{t}+T(\epsilon_1)} &\geq Q^*_{\bar{t}+T(\epsilon_1)-L}+\epsilon_1 +\sum_{i=1}^{T(\epsilon_1)-1} \left( Q^*_{\bar{t}+i-L}+\epsilon_1 - D_{\bar{t}+i} \right)\\
    & = \sum_{i=1}^{T(\epsilon_1)} Q^*_{\bar{t}+i-L} + T(\epsilon_1)\epsilon_1 - \sum_{i=1}^{T(\epsilon_1)-1} D_{\bar{t}+i}
\end{align*}
Assume that $T(\epsilon_1)>T(\epsilon_2)$. This means that for some $j<T(\epsilon_1)$ it must hold that 
\begin{equation*}
    D_{\bar{t}+j} \geq Q^*_{\bar{t}+j-L}+\epsilon_2 +\sum_{i=1}^{j-1} \left( Q^*_{\bar{t}+i-L}+\epsilon_2 - D_{\bar{t}+i} \right)
\end{equation*}
However, we know that $\epsilon_1<\epsilon_2$ and, therefore,
\begin{equation*}
    D_{\bar{t}+j} < Q^*_{\bar{t}+j-L}+\epsilon_1 +\sum_{i=1}^{j-1} \left( Q^*_{\bar{t}+i-L}+\epsilon_1 - D_{\bar{t}+i} \right) <  Q^*_{\bar{t}+j-L}+\epsilon_2 +\sum_{i=1}^{j-1} \left( Q^*_{\bar{t}+i-L}+\epsilon_2 - D_{\bar{t}+i} \right),
\end{equation*}
which leads to a contradiction.
Hence, for all $\epsilon_1$ and $\epsilon_2$, such that $\epsilon_1<\epsilon_2$, it must hold that $T(\epsilon_1)\leq T(\epsilon_2)$ and therefore the time between two stockouts under the $\epsilon$-perturbed policy is increasing in $\epsilon$. 
% \end{proof}
\hfill\Halmos
\endproof

\color{black}
Under the $\epsilon$-perturbed optimal policy, an additional $\epsilon$ units are ordered per period. This means that under the $\epsilon$-perturbed optimal policy there will always be more inventory than under the optimal policy. Hence, if a stockout is faced under the $\epsilon$-perturbed optimal policy, meaning that this larger inventory is depleted, all inventory under the optimal policy will be depleted. Therefore, if the $\epsilon$-perturbed optimal policy faces a stockout, then the optimal policy also faces a stockout. However, the reverse might not hold.
Therefore, let $Y_n(\epsilon)$ denote the number of stockouts under the optimal policy between two consecutive stockouts under the $\epsilon$-perturbed policy, i.e., $Y_n(\epsilon)$ denotes the number of stockouts under the optimal policy in the period $(S_{n-1}(\epsilon),S_n(\epsilon))$. 
In Lemma \ref{lem:Equal time between intermediate stockouts}, we show that, for a fixed value of $\epsilon$, the additional holding costs under the perturbed policy compared to the optimal policy are minimized when the intermediate stockouts under the optimal policy between stockouts of the perturbed policy occur at equal time intervals. In addition, in Lemma \ref{lem:coinciding_stockouts}, we show that the probability of such intermediate stockouts goes to zero as $\epsilon$ goes to zero.

\begin{lemma}\label{lem:Equal time between intermediate stockouts}
    Under continuous time, the additional holding costs under the $\epsilon$-perturbed policy compared to the optimal policy over the period $(S_{n-1}(\epsilon),S_n(\epsilon)]$ are minimized when the intermediate stockouts under the optimal policy between two consecutive stockouts under the $\epsilon$-perturbed policy occur at equally spread moments in time.
\end{lemma}

% \begin{proof}   
\proof{Proof.} 
In case that there are no intermediate stockouts, the difference in holding cost between the two policies equals $\sum_{j=1}^{T_n(\epsilon)-1}h\epsilon j$, where $T_n(\epsilon)$ denotes the time between the $(n-1)^{th}$ and $n^{th}$ stockout under the $\epsilon$-perturbed policy. 
In case of intermediate stockouts, some of these additional holding costs may be saved as part of the additional units that are ordered under the $\epsilon$-perturbed policy are used to satisfy demand. These cost savings are maximized when the additional units are exactly enough to satisfy the demand that would otherwise be lost.

Let $t_i$ denote the time of the $i^{th}$ intermediate stockout, $i=1,\ldots,y$. At time $t_1$, $t_1\epsilon$ additional units are available that can be used to satisfy demand. When these units are used, no holding costs occur in the remaining $T_n(\epsilon)-t_1$ periods, saving a total of $t_1\epsilon(T_n(\epsilon)-t_1)$ in holding costs compared to the case without intermediate stockouts. At the next intermediate stockout moment $t_2$, $(t_2-t_1)\epsilon$ additional units are available that can be used to satisfy demand, leading to $(t_2-t_1)\epsilon(T_n(\epsilon)-t_1)$ savings in additional holding costs, etc. This means that for $y$ intermediate stockouts, the savings in additional holding costs compared to the case without intermediate stockouts is given by 
\begin{equation*}
    t_1\epsilon(T_n-t_1) + \sum_{i=2}^y (t_i-t_{i-1})\epsilon(T_n(\epsilon)-t_i) = \epsilon\left(t_yT_n(\epsilon) - \sum_{i=1}^y t_i^2 + \sum_{i=2}^y t_{i-1}t_i\right).
\end{equation*}
We are interested in $t_1,\ldots,t_y$ that maximize these holding cost savings. To find these values for which the maximum is attained, we first derive the partial derivatives.
\begin{align*}
    \frac{\partial}{\partial t_1} \left[ t_yT_n(\epsilon) - \sum_{i=1}^y t_i^2 + \sum_{i=2}^y t_{i-1}t_i \right] &= -2t_1+t_2\\
    \frac{\partial}{\partial t_2} \left[ t_yT_n(\epsilon) - \sum_{i=1}^y t_i^2 + \sum_{i=2}^y t_{i-1}t_i \right] &= -2t_2+t_1+t_3\\
    \vdots & \\
    \frac{\partial}{\partial t_y} \left[ t_yT_n(\epsilon) - \sum_{i=1}^y t_i^2 + \sum_{i=2}^y t_{i-1}t_i \right] &= T_n(\epsilon) -2t_y+t_{y-1}
\end{align*}
Setting all partial derivatives equal to 0 and solving the system of equations gives the unique solution $t_1=\frac{1}{y+1}T_n(\epsilon),t_2=\frac{2}{y+1}T_n(\epsilon), \ldots, t_y=\frac{y}{y+1}T_n(\epsilon)$.
The Hessian matrix is given by:
\begin{equation*}
    H = \begin{bmatrix}
-2 & 1 & 0 & \ldots & 0 & 0 & 0 \\ 
1 & -2 & 1 & \ldots & 0 & 0 & 0 \\ 
\vdots & \vdots & \vdots & \ddots & \vdots & \vdots & \vdots \\ 
0 & 0 & 0 & \ldots & 1 & -2 & 1 \\
0 & 0 & 0 & \ldots & 0 & 1 & -2 
\end{bmatrix}
\end{equation*}
For any vector $\Vec{x}=[x_1 x_2 \ldots x_y]\in \mathbf{R}^y$, we find $\Vec{x}H\Vec{x}=-2 \sum_{i=1}^y x_i^2 + 2 \sum_{i=2}^y x_{i-1}x_i<0$. Therefore, the matrix $H$ is negative-definite and the given solution to the system of equations is indeed the (unique) maximum point.
Hence, the additional holding costs under the $\epsilon$-perturbed policy compared to the optimal policy are minimized when intermediate stockouts occur at equally spread moments in time.
% \end{proof}
\hfill\Halmos
\endproof

\begin{lemma}\label{lem:coinciding_stockouts}
    The probability that the stockout moments of the optimal policy and the $\epsilon$-perturbed policy do not coincide goes to zero as $\epsilon$ goes to zero.
\end{lemma}
% \begin{proof}
\proof{Proof.} 
We will proof this by induction. 
Let $I_t^*$ and $I_t^\epsilon$ denote the inventories at time $t$ under the optimal and the $\epsilon$-perturbed policy, respectively. Assume that a stockout is faced at time $t=0$ under both policies.

% \textbf{\textcolor{blue}{In het bewijs gebruiken we hier $S_1$, wat feitelijk afhankelijk is van epsilon. Daarom wilde ik het aanpassen naar $S_1(\epsilon)$ voor de eerste vergelijking, maar bij de tweede vergelijking moet het $S_1$ blijven omdat het voor het optimale policy is. In de inductie stap gaan we ervan uit dat $t=S_n$ een stockout moment is voor beide policies, waardoor ik niet meer helemaal zeker weet wat hier de beste notatie is. }} \textbf{\textcolor{red}{Kun je hier notatie $S_n(\epsilon)$ gebruiken? Dan hebben we $S_n(0)$ voor optimale policy. Lijkt hieronder te werken qua notatie.} \textcolor{blue}{Aangepast, moet nog checken of het zo klopt.}}

First, let $t=S_1(\epsilon)$ be the first moment after $t=0$ that a stockout is faced under the $\epsilon$-perturbed policy. Then the probability that there is no stockout under the optimal policy is denoted by 
\begin{align*}
    P\left(I_{S_1(\epsilon)}^* > 0 | I_{S_1(\epsilon)}^\epsilon=0 \right) &= P\left(\sum_{t=1}^{S_1(\epsilon)}Q_t^* > \sum_{t=1}^{S_1(\epsilon)} D_t | \sum_{t=1}^{S_1(\epsilon)}(Q_t^*+\epsilon) \leq \sum_{t=1}^{S_1(\epsilon)} D_t \right)=0.
\end{align*}
Alternatively, let $t=S_1(0)$ be the first moment after $t=0$ that a stockout is faced under the optimal policy. Then the probability that there is no stockout under the $\epsilon$-perturbed policy is denoted by 
\begin{align*}
    P\left(I_{S_1(0)}^\epsilon > 0 | I_{S_1(0)}^*=0 \right) &= P\left(\sum_{t=1}^{S_1(0)}(Q_t^*+\epsilon) > \sum_{t=1}^{S_1(0)} D_t | \sum_{t=1}^{S_1(0)}Q_t^* \leq \sum_{t=1}^{S_1(0)} D_t \right)\\
    &= P\left( S_1(0)\cdot\epsilon > \sum_{t=1}^{S_1(0)} (D_t - Q_t^*) | \sum_{t=1}^{S_1} (D_t - Q_t^*) \geq 0 \right)\\
    &< P\left( S_1(0)\cdot\epsilon > 0 \right).
\end{align*}
This probability is monotonically decreasing to 0 as $\epsilon\rightarrow0$. Hence, the probability that there is no stockout at time $S_1(0)$ under the $\epsilon$-perturbed policy when there is a stockout under the optimal policy goes to 0 as $\epsilon\rightarrow0$.

Now, assume that there is a stockout under both the optimal policy and the $\epsilon$-perturbed policy at time $t=S_n(0)=S_n(\epsilon)$. %\textbf{\textcolor{blue}{Already claiming too much by saying t is the n-th stockout moment under both policies?}}
First, let $t=S_{n+1}(\epsilon)$ be the next stockout moment under the $\epsilon$-perturbed policy. Then the probability that there is no stockout under the optimal policy is denoted by 
\begin{align*}
    P\left(I_{S_{n+1}(\epsilon)}^* > 0 | I_{S_{n+1}(\epsilon)}^\epsilon=0 \right) &= P\left(\sum_{t=S_n(\epsilon)+1}^{S_{n+1}(\epsilon)}Q_t^* > \sum_{t=S_n(\epsilon)+1}^{S_{n+1}(\epsilon)} D_t | \sum_{t=S_n(\epsilon)+1}^{S_{n+1}(\epsilon)}(Q_t^*+\epsilon) \leq \sum_{t=S_n(\epsilon)+1}^{S_{n+1}(\epsilon)} D_t \right)= 0.
\end{align*}
Now assume that $t=S_{n+1}(0)$ is the next stockout moment after $S_n$ under the optimal policy. Then,
\begin{align*}
    P\left(I_{S_{n+1}(0)}^\epsilon > 0 | I_{S_{n+1}(0)}^*=0 \right) &= P\left(\sum_{t=S_n(0)+1}^{S_{n+1}(0)}(Q_t^*+\epsilon) > \sum_{t=S_n(0)+1}^{S_{n+1}(0)} D_t | \sum_{t=S_n(0)+1}^{S_{n+1}(0)}Q_t^* \leq \sum_{t=S_n(0)+1}^{S_{n+1}(0)} D_t \right)\\
    &= P\left( (S_{n+1}(0)-S_n(0))\cdot\epsilon > \sum_{t=S_n(0)+1}^{S_{n+1}(0)} (D_t - Q_t^*) | \sum_{t=S_n(0)+1}^{S_{n+1}(0)} (D_t - Q_t^*) \geq 0 \right)\\
    &< P\left( (S_{n+1}(0)-S_n(0))\cdot\epsilon > 0 \right).
\end{align*}
Again, this probability is monotonically decreasing to 0 as $\epsilon\rightarrow0$. Therefore, if $S_n(\epsilon)=S_n(0)$ is a stockout moment under both policies, the probability that there is no stockout under the $\epsilon$-perturbed policy the next time there is a stockout under the optimal policy goes to zero as $\epsilon\rightarrow0$. 

Hence, using induction, it can be concluded that the probability that a stockout is faced under the optimal policy but not under the $\epsilon$-perturbed policy goes monotonically to zero as $\epsilon\rightarrow0$.
Consequently, as $\epsilon$ goes to zero, the probability that the stockout moments under the optimal policy and the $\epsilon$-perturbed policy do not coincide goes monotonically to zero.
% \end{proof}
\hfill\Halmos
\endproof

With the properties of the $\epsilon$-perturbed optimal policy, we are able to proof a property of the optimal policy without knowing the optimal policy itself. This property can be used to test if the optimal policy within a policy class can be globally optimal. In particular, if the optimal policy within a class of policies does not satisfy the property, then the class of policies does not contain the globally optimal policy. The property follows by deriving an explicit expression for the difference in long-run average cost between the $\epsilon$-perturbed policy and the optimal policy and using the KKT condition as we let $\epsilon$ approach zero.   

\color{black}
\begin{theorem}{\textbf{Optimality equation for lost-sales systems}}\label{lem:Optimality equation for lost-sales systems}
\newline
    Under the optimal policy the following equation holds:
    \begin{equation}
\label{eq1}
\frac{E[T^2]}{E[T]} = {\frac{2p+h}{h}}\\
\end{equation}
\end{theorem}

% \begin{proof}
\proof{Proof.} 
Define
\begin{center}
    $C^*$ := long-run average cost under the optimal policy, and
\\
$C_\epsilon$ := long-run average cost under the $\epsilon$-perturbed policy.
\end{center}

Let us consider the cost under both policies over the interval $(S_{n-1}(\epsilon),S_n(\epsilon)]$. 
If $Y_n(\epsilon)=0$, i.e., the stockout moments under both policies coincide, at time $S_{n-1}(\epsilon)$ the inventory immediately after arrival of the replenishment order equals $Q^*_{S_{n-1}(\epsilon)-L}$ under the optimal policy, and equals $Q^*_{S_{n-1}(\epsilon)-L}+\epsilon$ under the $\epsilon$-perturbed optimal policy. Until time $S_n(\epsilon)$ each order under the $\epsilon$-perturbed optimal policy adds another $\epsilon$ to the inventory in comparison to the optimal policy. At time $S_n(\epsilon)$ a stockout occurs. Since we added $T_n(\epsilon)$ times $\epsilon$ to the inventory under the optimal policy, the amount lost under the $\epsilon$-perturbed policy is $T_n(\epsilon)\epsilon$ less than under the optimal policy. Thus, the cost difference over the interval $(S_{n-1}(\epsilon),S_n(\epsilon)]$ between the two policies equals
% \textbf{\textcolor{blue}{We zijn volgens mij niet helemaal consequent met het gebruik van $k$. We gebruiken het deels om opeenvolgende periodes aan te geven, zoals in de eerste term hieronder waar ik het nu aangepast heb naar $j$, en (groten)deels om de stockout intervallen van het epsilon policy te tellen (zoals gedefinieerd helemaal aan het begin). Ik heb nu in het eerste geval k aangepast naar j, maar dit is niet optimaal want we gebruiken j ook al in het bewijs van Lemma 1. Misschien kunnen we gewoon t gebruiken? }} \textbf{\textcolor{red}{Er zijn vast andere indices. We reserveren t voor de beslistijdstippen. Dus bij opeenvolgende perioden (dus met 1 ertussen) kun je dus t gebruiken. t is hiermee een soort basisindex. Als je dan sommeert over opeenvolgende perioden, gebruik je een andere index, omdat vrijwel overal t als index in staat. } \textcolor{blue}{Ik heb de index nu aangepast naar $m$, deze gebruiken we volgens mij niet meer ergens anders.}}
\begin{equation*}
\sum_{m=1}^{T_n(\epsilon)-1}{h\epsilon m} - pT_n(\epsilon)\epsilon = h\epsilon\frac{T_n(\epsilon)(T_n(\epsilon)-1)}{2}-p\epsilon T_n(\epsilon).\\
\end{equation*}

If $Y_n(\epsilon)>0$, i.e., there are additional stockouts under the optimal policy between two consecutive stockouts under the $\epsilon$-perturbed policy, the cost difference decreases. The amount lost under the $\epsilon$-perturbed policy is still $T_n(\epsilon)\epsilon$ less than under the optimal policy, as the total difference in orders placed over the period $(S_{n-1}(\epsilon),S_n(\epsilon)]$ remains the same. However, the additional holding costs under the $\epsilon$-perturbed policy compared to the optimal policy are smaller than $\sum_{m=1}^{T_n(\epsilon)-1}{h\epsilon m}$, as some of these additional orders are already used to satisfy demand. 

For known combinations of $T_n(\epsilon)$ and $Y_n(\epsilon)$ an upper bound on the additional holding costs that can be saved compared to the no-intermediate-stockout scenario can be computed. 
Under continuous time, the cost savings would be maximized when the intermediate stockouts are evenly spread out over the time between two consecutive stockouts of the $\epsilon$-perturbed policy (see Lemma 1), i.e., for $Y_n(\epsilon)=y$ the intermediate stockouts occur at time $\frac{1}{y+1}T_n(\epsilon),\frac{2}{y+1}T_n(\epsilon),\ldots,\frac{y}{y+1}T_n(\epsilon)$. Note that the division of the time period $T_n(\epsilon)$ into equal sub-intervals may lead to non-integer time points. Even though we consider discrete time, this can be disregarded, as rounding to integer time points can only reduce the potential cost savings below this upper bound.
Therefore, by considering all $y$ intermediate stockouts, an upper bound on the total additional holding cost savings in the intermediate stockout scenario compared to the no-intermediate stockout scenario can be written as
\begin{multline*}
    \frac{T_n(\epsilon)}{y+1}\sum_{i=1}^{y}\left(T_n(\epsilon)-i\frac{T_n(\epsilon)}{y+1}\right)\epsilon h=
    \frac{T_n(\epsilon)}{y+1}\left(yT_n(\epsilon)-\frac{T_n(\epsilon)}{y+1}\frac{y(y+1)}{2}\right)\epsilon h \\=
    \frac{T_n(\epsilon)}{y+1}\left(yT_n(\epsilon)-\frac{y}{2}T_n(\epsilon)\right)\epsilon h =
    \frac{y}{y+1}\frac{T_n(\epsilon)^2}{2}\epsilon h.
\end{multline*}

We know that the number of additional stockouts is bounded by the number of periods between the consecutive stockouts under the $\epsilon$-perturbed policy. Therefore, in this case, the cost difference between the optimal and $\epsilon$-perturbed policy over the interval $(S_{n-1}(\epsilon),S_n(\epsilon)]$ between the two policies equals
\begin{align*}
&\sum_{m=1}^{T_n(\epsilon)-1}{h\epsilon m} - \sum_{y=1}^{T_n(\epsilon)-1}\left[\frac{y}{y+1}\frac{T_n(\epsilon)^2}{2}\epsilon h P\{Y_n(\epsilon)=y\} \right]  - pT_n(\epsilon)\epsilon \\
&= h\epsilon\frac{T_n(\epsilon)(T_n(\epsilon)-1)}{2}-p\epsilon T_n(\epsilon) - \epsilon h \sum_{y=1}^{T_n(\epsilon)-1}\left[ P\{Y_n(\epsilon)=y\} \frac{y}{y+1}\frac{T_n(\epsilon)^2}{2} \right].
\end{align*}

We are interested in the total difference in holding and penalty costs under the optimal policy and the $\epsilon$-perturbed policy. 
\color{black}
It follows that in case all stockouts coincide,
\begin{equation*}
\begin{aligned}
C_\epsilon - C^* & = \lim_{n\to\infty} \frac {\sum_{k=1}^n \left[ h\epsilon\frac{T_k(\epsilon)(T_k(\epsilon)-1)}{2}-p\epsilon T_k(\epsilon) \right]}{\sum_{k=1}^n T_k(\epsilon)}\\
& =\frac{\lim_{n\to\infty}\frac{1}{n} \sum_{k=1}^n{\epsilon\left(\frac{h}{2}T_k(\epsilon)(T_k(\epsilon)-1)-pT_k(\epsilon)\right)}}{\lim_{n\to\infty}\frac{1}{n}\sum_{k=1}^n T_k(\epsilon)}\\
& =\epsilon \frac{\frac{h}{2} E[T(T-1)]-pE[T]}{E[T]}.\\
\end{aligned}    
\end{equation*}

Alternatively, if there are additional stockouts faced under the optimal policy that are not faced under the $\epsilon$-perturbed policy, we find
\begin{equation*}
\begin{aligned}
C_\epsilon - C^*  \geq & \lim_{n\to\infty} \frac {
\sum_{k=1}^n \left[ h\epsilon\frac{T_k(\epsilon)(T_k(\epsilon)-1)}{2}-p\epsilon T_k(\epsilon) - \sum_{y=1}^{T_k(\epsilon)-1}\left[ \frac{y}{y+1}\frac{T_k(\epsilon)^2}{2}\epsilon h P\{Y_k(\epsilon)=y\} \right] \right]
}{\sum_{k=1}^n T_k}\\
=& \lim_{n\to\infty} \frac { \frac{1}{n}\sum_{k=1}^n \left[ h\epsilon\frac{T_k(\epsilon)(T_k(\epsilon)-1)}{2}-p\epsilon T_k(\epsilon) \right] }{\frac{1}{n}\sum_{k=1}^n T_k(\epsilon)}  - 
\lim_{n\to\infty} \frac {\frac{1}{n}\sum_{k=1}^n\sum_{y=1}^{T_k(\epsilon)-1}\left[\frac{y}{y+1}\frac{T_k(\epsilon)^2}{2}\epsilon h P\{Y_k(\epsilon)=y\} \right]}{\frac{1}{n}\sum_{k=1}^n T_k(\epsilon)}\\
 = &\frac{\lim_{n\to\infty}\frac{1}{n} \sum_{k=1}^n{\epsilon(\frac{h}{2}T_k(\epsilon)(T_k(\epsilon)-1)-pT_k(\epsilon))}}{\lim_{n\to\infty}\frac{1}{n}\sum_{k=1}^n T_k(\epsilon)} - 
 \frac {\lim_{n\to\infty}\frac{1}{n}\sum_{k=1}^n\sum_{y=1}^{T_k(\epsilon)-1}\left[ \frac{y}{y+1}\frac{T_k(\epsilon)^2}{2}\epsilon h P\{Y_k(\epsilon)=y\} \right]}{\lim_{n\to\infty}\frac{1}{n}\sum_{k=1}^n T_k(\epsilon)}\\
=& \epsilon \frac{\frac{h}{2} E[T(T-1)]-pE[T]}{E[T]} - 
\frac {\lim_{n\to\infty}\frac{1}{n}\sum_{k=1}^n\sum_{y=1}^{T_k(\epsilon)-1}\left[\frac{y}{y+1}\frac{T_k(\epsilon)^2}{2}\epsilon h P\{Y_k(\epsilon)=y\} \right]}{\lim_{n\to\infty}\frac{1}{n}\sum_{k=1}^n T_k(\epsilon)}
\end{aligned}    
\end{equation*}

In general, considering that there may or may not be intermediate stockouts, we thus find that 
\begin{multline*}
\epsilon \frac{\frac{h}{2} E[T(T-1)]-pE[T]}{E[T]} - 
\frac {\lim_{n\to\infty}\frac{1}{n}\sum_{k=1}^n\sum_{y=1}^{T_k(\epsilon)-1}\left[ \frac{y}{y+1}\frac{T_k(\epsilon)^2}{2}\epsilon h P\{Y_k(\epsilon)=y\} \right]}{\lim_{n\to\infty}\frac{1}{n}\sum_{k=1}^n T_k(\epsilon)} \\ \leq C_{\epsilon} - C^* \leq\epsilon \frac{\frac{h}{2} E[T(T-1)]-pE[T]}{E[T]}   ,
\end{multline*}
and consequently,
\begin{multline*}
\frac{\frac{h}{2} E[T(T-1)]-pE[T]}{E[T]} - \lim_{\epsilon\to 0} 
\frac {\lim_{n\to\infty}\frac{1}{n}\sum_{k=1}^n\sum_{y=1}^{T_k(\epsilon)-1}\left[ \frac{y}{y+1}\frac{T_k(\epsilon)^2}{2} h P\{Y_k(\epsilon)=y\} \right]}{\lim_{n\to\infty}\frac{1}{n}\sum_{k=1}^n T_k(\epsilon)} \\ \leq \lim_{\epsilon\to 0}\frac{C_{\epsilon} - C^*}{\epsilon} \leq \frac{\frac{h}{2} E[T(T-1)]-pE[T]}{E[T]}.  
\end{multline*}

In order to show that $\lim_{\epsilon\to 0} \frac{C_{\epsilon} - C^*}{\epsilon} = \frac{\frac{h}{2} E[T(T-1)]-pE[T]}{E[T]}$, we take a closer look at 
\begin{align*}
    &\lim_{\epsilon\to 0} 
\frac {\lim_{n\to\infty}\frac{1}{n}\sum_{k=1}^n\sum_{y=1}^{T_k(\epsilon)-1}\left[ \frac{y}{y+1}\frac{T_k(\epsilon)^2}{2} h P\{Y_k(\epsilon)=y\} \right]}{\lim_{n\to\infty}\frac{1}{n}\sum_{k=1}^n T_k(\epsilon)}\\
&= 
h \frac {\lim_{\epsilon\to 0}\lim_{n\to\infty}\frac{1}{n}\sum_{k=1}^n\sum_{y=1}^{T_k(\epsilon)-1}\left[ \frac{y}{y+1}\frac{T_k(\epsilon)^2}{2} P\{Y_k(\epsilon)=y\}\right]}{\lim_{\epsilon\to 0}\lim_{n\to\infty}\frac{1}{n}\sum_{k=1}^n T_k(\epsilon)}.
\end{align*}
Looking at numerator, we need to show that:
\begin{align*}
    \lim_{\epsilon\to 0}\lim_{n\to\infty}\frac{1}{n}\sum_{k=1}^n \frac{T_k(\epsilon)^2}{2}\sum_{y=1}^{T_k(\epsilon)-1}\left[ \frac{y}{y+1} P\{Y_k(\epsilon)=y\} \right] &= 0. 
\end{align*}

% \color{red}
% \textbf{\textcolor{blue}{In rood volgt het deel dat we bij ISIR besproken hebben. Volgens mij vult dat inderdaad het gat dat we nog hadden. Het kan zijn dat het nog iets compacter opgeschreven kan worden dan ik nu heb gedaan.}}
We know that
\begin{align*}
    \lim_{\epsilon\to 0}\lim_{n\to\infty}\frac{1}{n}\sum_{k=1}^n \frac{T_k(\epsilon)^2}{2}\sum_{y=1}^{T_k(\epsilon)-1}\left[ \frac{y}{y+1} P\{Y_k(\epsilon)=y\} \right] &\leq 
    \lim_{\epsilon\to 0}\lim_{n\to\infty}\frac{1}{n}\sum_{k=1}^n \frac{T_k(\epsilon)^2}{2}\sum_{y=1}^{T_k(\epsilon)-1}\left[ P\{Y_k(\epsilon)=y\} \right] \\
    &= 
    \lim_{\epsilon\to 0}\lim_{n\to\infty}\frac{1}{n}\sum_{k=1}^n \frac{T_k(\epsilon)^2}{2} P\{Y_k(\epsilon)\geq1\}\\
    &=
     \lim_{\epsilon\to 0} \left[ P\{Y(\epsilon)\geq1\} \lim_{n\to\infty}\frac{1}{n}\sum_{k=1}^n \frac{T_k(\epsilon)^2}{2} \right],
\end{align*}
where the last equality holds because the probability of intermediate stockouts occurring is independent of the time interval $k$.
In Lemma \ref{lem:time_between_stockouts_increasing_in_eps} we have shown that the time between stockouts under the $\epsilon-$perturbed policy is increasing in $\epsilon$. Therefore, we can bound $T_k^2(\epsilon)$ by $T_k^2(\epsilon_{max})$ for each interval $k$, where $\epsilon_{max}$ denotes the maximum value $\epsilon$ can take. 
To assure that the system remains stable, we set $\epsilon_{max}=\frac{\mathbb{E}[D]-\mathbb{E}[Q^*]}{2}$. 
% \textcolor{blue}{Moeten we dit hier introduceren of aan het begin van de sectie, waar we het epsilon-perturbed policy met epsilon$>$0 introduceren?}}\textbf{\textcolor{red}{Alles wat hiervoor stond is waar zonder de restrictie. Door het hier te zetten, zie je waarom dit nu nodig is: om te zorgen dat we over stabiele sytemen spreken.}}\\
Therefore, it holds that 
\begin{align*}
     \lim_{\epsilon\to 0} \left[ P\{Y(\epsilon)\geq1\} \lim_{n\to\infty}\frac{1}{n}\sum_{k=1}^n \frac{T_k(\epsilon)^2}{2} \right] &\leq  
    \lim_{\epsilon\to 0} \left[ P\{Y(\epsilon)\geq1\} \lim_{n\to\infty}\frac{1}{n}\sum_{k=1}^n \frac{T_k(\epsilon_{max})^2}{2} \right].
\end{align*}
Since $\lim_{n\to\infty}\frac{1}{n}\sum_{k=1}^n \frac{T_k(\epsilon_{max})^2}{2}=\mathbb{E}\left[ \frac{1}{2} T(\epsilon_{max})^2 \right]$ and this limit is finite,  
\begin{align*}
     \lim_{\epsilon\to 0} \left[ P\{Y(\epsilon)\geq1\} \lim_{n\to\infty}\frac{1}{n}\sum_{k=1}^n \frac{T_k(\epsilon_{max})^2}{2} \right] &= 
    \lim_{\epsilon\to 0} \left[ P\{Y(\epsilon)\geq1\} \mathbb{E}\left[ \frac{1}{2} T(\epsilon_{max})^2 \right] \right]\\
    &= 
     \lim_{\epsilon\to 0} \left[ P\{Y(\epsilon)\geq1\}\right] \cdot \mathbb{E}\left[ \frac{1}{2} T(\epsilon_{max})^2 \right] \\
     & = 0 \cdot \mathbb{E}\left[ \frac{1}{2} T(\epsilon_{max})^2 \right] = 0,
\end{align*}
as we have shown in Lemma \ref{lem:coinciding_stockouts} that the probability that the stockout moments of the optimal policy and the $\epsilon-$perturbed policy do not coincide goes to zero as $\epsilon$ goes to zero.
Combining above steps with the fact that all terms are non-negative, this means that 
\begin{align*}
    \lim_{\epsilon\to 0}\lim_{n\to\infty}\frac{1}{n}\sum_{k=1}^n \frac{T_k(\epsilon)^2}{2}\sum_{y=1}^{T_k(\epsilon)-1}\left[ \frac{y}{y+1} P\{Y_k(\epsilon)=y\} \right] &= 0.
\end{align*}

Hence, we find
\begin{equation*}
\begin{aligned}
\frac{\frac{h}{2} E[T(T-1)]-pE[T]}{E[T]} - 0h  \leq \lim_{\epsilon\to 0}\frac{C_{\epsilon} - C^*}{\epsilon} \leq \frac{\frac{h}{2} E[T(T-1)]-pE[T]}{E[T]}\\
\end{aligned}    
\end{equation*}
\color{black}
From this, we derive 
\begin{equation*}
\begin{aligned}
\lim_{\epsilon\to 0}\frac{C_{\epsilon} - C^*}{\epsilon} = \frac{\frac{h}{2} E[T(T-1)]-pE[T]}{E[T]}= \frac{h}{2} \frac{E[T^2]}{E[T]} - \frac{h}{2} -p.
\end{aligned}    
\end{equation*}
As $C^*$ is the minimal cost under the optimal policy, we also have
\begin{equation*}
\begin{aligned}
\lim_{\epsilon\to 0}\frac{C_{\epsilon} - C^*}{\epsilon} = 0.
\end{aligned}    
\end{equation*}
Together, this implies that
\begin{equation*}
\begin{aligned}
\frac{h}{2} \frac{E[T^2]}{E[T]} = \frac{h}{2} + p,
\end{aligned}    
\end{equation*}
This proves the theorem.
% \end{proof}
\hfill\Halmos
\endproof

Note that $\frac{E[T^2]}{E[T]}$ is equal to the average time between the stockouts before and after some time $t$ for $t\to\infty$, if the times between stockouts would constitute a renewal process.  It immediately follows from the optimality equation that $\frac{E[T^2]}{E[T]}$ is increasing in the penalty cost $p$. With increasing penalty cost, the non-stockout probability under the optimal policy increases. This suggests that $\frac{E[T^2]}{E[T]}$ is increasing in the non-stockout probability. Though we do not provide a formal proof of these statements, we have exploited them in our numerical experiments to find close-to optimal parameters for the fixed non-stockout probability ($FP_3$) policy. In the next section we discuss both $FP_3$-policy and $PIL$-policy in detail, deriving both exact expressions and computationally efficient approximations. 

%% \[\lim_{x\to\infty} f(x) \]

\section{Myopic policies for the lost-sales problem}
\label{sec4}
In \cite{Dons1996} it is shown that for the lost-sales model the $FP_3$-policy is superior to a base stock policy: the myopic policy needs a lower average inventory and has a lower replenishment order volatility. For an inventory model with backordering and linear holding and penalty costs, \cite{Kok2018} shows that under mild conditions on the policy structure, the optimal policy within a class of policies satisfies the Newsvendor fractile, implying that the non-stockout probability at an arbitrary point in time must satisfy the Newsvendor fractile $\frac{p}{p+h}$. For the periodic review backlog model without fixed costs and linear holding and penalty costs, the optimal policy is a myopic policy that aims to meet this Newsvendor fractile every period, which is in fact equivalent to a base stock policy. It should be noted that the cost structure for this classical backlog model is different from the cost structure for the lost-sales model in this paper. In the lost-sales model the penalty costs are linear in the number of items lost at the end of a period, whereas the penalty costs in the backlog model are linear in the number of items backlogged at the end of a period. While in the lost-sales model an item can incur penalty costs only once, in the backlog model the same item may be subject to a penalty cost over a number of periods. Still, the above discussion suggests that it should be worthwhile to study the cost performance of the myopic fixed non-stockout probability policy, where the replenishment order released ensures that at the moment of release, the non-stockout probability at the end of the period at the start of which this order replenishes inventory equals a fixed target value $P_3^*$.

\subsection{Fixed non-stockout probability policy}
\label{subsec41}

The non-stockout probability at time $t+L+1$ can be written as L subsequent Lindley equations 
\begin{equation*}
    P_{3,t+L}=P\{\tilde{I}_{t+L}>0\}.
\end{equation*}
It is easy to see that
\begin{equation*}
    \tilde{I}_t=(I_t-D_t)^+
\end{equation*}
\begin{equation}
\label{eqbal}
    \tilde{I}_{t+s}=\tilde{I}_{t+s-1}+Q_{t-(L-s)}-D_{t+s}, s=1,...,L.
\end{equation}
Then we write the expression for $P_{3,t+L}$ in terms of $I_t, D_{t+s}, s=0,...,L$ and $Q_{t-L+s}, s=1,...,L$ as follows.
\begin{equation}
    P_{3,t+L}=P\{(...((I_t-D_t)^+ + Q_{t-L+1}-D_{t+1})^+ + ...)^++Q_t-D_{t+L} > 0\}
\end{equation}
This non-stockout probability is determined at time $t$, when ordering $Q_t$. Under the fixed non-stockout probability policy with target $P_3^*$, $\tilde{Q_t}$ is determined by
\begin{equation}\label{eq: non-stockout probability}
    P\{(...((I_t-D_t)^+ + Q_{t-L+1}-D_{t+1})^+ + ...)^++\tilde{Q_t}-D_{t+L} > 0\} = P_3^*.
\end{equation}

It is easy to prove that, when starting with no outstanding orders and zero inventory, under continuous demand it is possible for all $t\geq0$ to find $\tilde{Q_t}$ as a solution of equation \eqref{eq: non-stockout probability}. In \cite{Dons1996} a numerical procedure is proposed to approximately solve equation \eqref{eq: non-stockout probability}, using the moment-iteration procedure for the G/G/1 queue proposed by \cite{Kok1989}.  We denote this procedure as the forward recursion. As it turns out, an alternative procedure can be applied to solve equation \eqref{eq: non-stockout probability}, along the same lines as proposed in \cite{Kok2003} for computing ruin probabilities. Firstly, this yields the following identities.
\begin{equation*}
\begin{aligned}
    P_{3,t+L}&=P\{(...((I_t-D_t)^+ + Q_{t-L+1}-D_{t+1})^+ + ...)^++Q_t-D_{t+L} > 0\}\\
&=1-P\{((...((I_t-D_t)^+ + Q_{t-L+1}-D_{t+1})^+ + ...)^++Q_t-D_{t+L})^+ = 0\}\\
&=1-P\{D_{t+L}>Q_t,D_{t+L+1}+D_{t+L}>Q_t+Q_{t-1},...,\sum_{s=0}^{L}D_{t+s}>I_t+\sum_{s=0}^{L-1}Q_{t-s}\}
\end{aligned}    
\end{equation*}
Though the last equation shows that $P_{3,t+L}$ is determined by $L+1$ mutually dependent events, we can derive a set of recursive equations, introducing conditional random variables, that can be used for an efficient and accurate approximation of $P_{3,t+L}$. Towards that end we define,
\begin{equation*}
\begin{aligned}
\xi_0 &:= 0\\
\xi_n &:= \sum_{m=1}^nQ_{t+1-m}, n=1,...,L\\
\xi_{L+1}&:=I_t+\xi_L\\
X_n&:=D_{t+L+1-n}, n=1,...,L+1\\
\end{aligned}    
\end{equation*}
Then it follows that
\begin{equation}
\label{eqP3}
\begin{aligned}
    P_{3,t+L}&=1-P\{\sum_{m=1}^nX_m>\xi_n, n=1,...,L+1\}.
\end{aligned}    
\end{equation}
With this we define the random variables $Y_n, n=1,...,L+1$, to satisfy
\begin{equation}
\label{eqYn}
\begin{aligned}
P\{Y_n>x\}= \frac{P\{\sum_{k=1}^mX_k>\xi_m, m=1,...,n-1,\sum_{k=1}^nX_k>x\}}{P\{\sum_{k=1}^mX_k>\xi_m, m=1,...,n-1\}}.
\end{aligned}    
\end{equation}
Note that
\begin{equation*}
\begin{aligned}
P\{Y_1>x\}= P\{X_1>x\}.
\end{aligned}    
\end{equation*}
Similar as in \cite{Kok2003} we can derive from the above definition of $Y_n, {n=1,...,L+1,}$ that
\begin{equation*}
\begin{aligned}
Y_n=X_n+(Y_{n-1}|Y_{n-1}>\xi_{n-1}), n=1,...,L+1,
\end{aligned}    
\end{equation*}
which implies that
\begin{equation*}
\begin{aligned}
P\{Y_n>\xi_{n-1}\}=1, n=1,...,L+1.
\end{aligned}    
\end{equation*}
We can define non-negative random variables $Z_n, n=1,...,L+1$,
\begin{equation}
\label{eqZn}
Z_n=Y_n- \xi_{n-1}, n=1,...,L+1.   
\end{equation}
Then we find the following recursion for the conditional random variables $Z_n, n=1,...,L+1$,
\begin{equation*}
\begin{aligned}
Z_1 &= X_1,\\
Z_n &= X_n+(Z_{n-1}-(\xi_{n-1}-\xi_{n-2}|Z_{n-1}>\xi_{n-1}-\xi_{n-2}), n=2,...,L+1.
\end{aligned}    
\end{equation*}
This exact recursion is the basis for an efficient computation of an approximation of $P_{3,t+L}$, which is described in subsection \ref{algp3}. 

\subsection{Asymptotic behaviour of the $FP_3$-policy}
\label{subsec42}
In the next section we show the accuracy of the approximation for $P_{3,t+L}$ and discuss the performance of the $FP_3$ policy against that of asymptotically optimal policies. In this subsection we discuss the asymptotic behaviour of the $FP_3$-policy for high target values of $P_{3,t+L}$ and long lead times. Let us first consider the situation of high target values of $P_{3,t+L}$. In order to maintain such high service levels, we need to maintain high levels of $I_t$. We can see the implication of high levels of $I_t$ on the decision variable $Q_t$ by rewriting equation \eqref{eqP3} as follows,
\begin{equation*}
\begin{aligned}
    P_{3,t+L}&=1-P\{\sum_{m=1}^nX_m>\xi_n, n=1,...,L+1\}\\
    &=P\{\cup_{n=1}^{L+1}\{\sum_{m=1}^nX_m\leq\xi_n\}\}\\
    &=P\{\sum_{m=1}^{L+1}X_m\leq\xi_{L+1}\}+P\{\sum_{m=1}^{L+1}X_m>\xi_{L+1},\cup_{n=1}^{L}\{\sum_{m=1}^nX_m\leq\xi_n\}\}.
\end{aligned}
\end{equation*}

As the penalty cost $p\to\infty$, $P_3^* \to 1$ and thereby $I_t \to \infty$. This implies that $\xi_{L+1} \to \infty$. Then it follows that
\begin{equation*}
\begin{aligned}
\lim_{p\to\infty} \frac {P_{3,t+L}}{P\{\sum_{m=1}^{L+1}X_m\leq\xi_{L+1}\}}=1.
\end{aligned}
\end{equation*}
Hence for large values of $p$ we can solve for $Q_t$ by solving
\begin{equation*}
\begin{aligned}
{P\{\sum_{s=0}^LD_{t+s} \leq I_t+\sum_{s=0}^{L-1}Q_{t-s}\}}=P_3^*.
\end{aligned}
\end{equation*}
Solving this equation is equivalent to operating under a base stock policy. Thus, for high penalty costs the $FP_3$-policy behaves like a base stock policy.

Now let us consider the case for $L\to\infty$. Assuming that under the $FP_3$-policy the system eventually reaches stationarity, the order quantity $Q_t$ arrives after an infinite amount of time and therefore finds the inventory in its stationary state $I$. Thus $Q_t$ follows from the following equation,
\begin{equation*}
\begin{aligned}
P\{I+Q_t-D \ge0\}=P_3^*.
\end{aligned}
\end{equation*}
Then it follows that every $Q_t$ satisfies the same equation and thus for $L\to\infty$, the $FP_3$-policy behaves as a constant order policy. Thus the $FP_3$-policy is asymptotically optimal under the asymptotic regimes studied in \cite{huh2009asymptotic} and \cite{goldberg2016asymptotic}.

\subsection{Projected inventory level policy}
\label{subsec43}
A myopic policy related to the fixed-fill-rate and fixed-non-stockout-probability policies is the fixed projected-inventory-level ($PIL$) policy proposed by \cite{van2024projected}. Under this policy, each period the replenishment order quantity is determined such that a target expected inventory at the start of the period, immediately after the order is added to inventory, is met. It is shown in \cite{van2024projected} that the $PIL$ policy is asymptotically optimal for long lead times and high penalty costs. In section \ref{sec6} we compare the $FP_3$ policy with the $PIL$ policy. For that purpose we develop an exact expression for the expected inventory $E[\tilde{I}_{t+L-1}]$. Note that
\begin{equation*}
\tilde{I}_{t+L-1} := (I_{t+L-1}-D_{t+L-1})^+.
\end{equation*}
Then it follows that
\begin{equation*}
I_{t+L}=\tilde{I}_{t+L-1}+Q_t.\\
\end{equation*}
Basic probability theory yields
\begin{equation}\label{eq: exp inv}
E[\tilde{I}_{t+L-1}] := \int_0^\infty(1-P\{\tilde{I}_{t+L-1} \leq x \})dx.
\end{equation}
Similar to the derivation of the expression for $P_{3,t+L}$, we find an expression for $P\{\tilde{I}_{t+L-1} \leq x \}$,
\begin{equation*}
\begin{aligned}
P\{\tilde{I}_{t+L-1} \leq x \} &= P\{(...((I_t-D_t)^+ + Q_{t-L+1}-D_{t+1})^+ + ...)^++Q_{t-1}-D_{t+L-1})^+ \leq x\}\\
&= P\{(...((I_t-D_t)^+ + Q_{t-L+1}-D_{t+1})^+ + ...)^++Q_{t-1}-x-D_{t+L-1})^+ =0\}
\end{aligned}
\end{equation*}

We define $\tilde{X}_n$ and $\tilde{\xi}_n, n=1,...,L$,
\begin{equation*}
\begin{aligned}
\tilde{\xi}_0 &:= 0\\
\tilde{\xi}_n &:= \sum_{m=1}^nQ_{t-m}, n=1,...,L-1\\
\tilde{\xi}_{L}&:=I_t+\tilde{\xi}_{L-1}\\
\tilde{X}_n&:=D_{t+L-n}, n=1,...,L\\
\end{aligned}    
\end{equation*}
Similar to the derivation of equation \eqref{eqP3} we find
\begin{equation*}
P\{\tilde{I}_{t+L-1} \leq x \} = P\{\sum_{m=1}^{n}\tilde{X}_m \geq \tilde{\xi}_n-x,n=1,...,L\}.\\
\end{equation*}
Substitution of this equation into equation \eqref{eq: exp inv} and using $\tilde{\xi}_0=0$ yields
\begin{equation*}
\begin{aligned}
E[\tilde{I}_{t+L-1}] &= \int_0^\infty(1-P\{\sum_{m=1}^{n}\tilde{X}_m \geq \tilde{\xi}_n-x,n=1,...,L\})dx\\
&= \sum_{k=1}^{L}{\int_{\tilde{\xi}_{k-1}}^{\tilde{\xi}_k}(1-P\{\sum_{m=1}^{n}\tilde{X}_m \geq \tilde{\xi}_n-x,n=1,...,L\})dx}.   
\end{aligned}
\end{equation*}
As $\tilde{\xi}_1 \leq \tilde{\xi_2} \leq ... \leq\tilde{\xi}_{L}$, we find that
\begin{equation*}
\begin{aligned}
\tilde{\xi}_n-x \leq 0, n=0,...,k-1, x \in (\tilde{\xi}_{k-1},\tilde{\xi}_k), k=1,...,L.
\end{aligned}\end{equation*}
With this we can write $E[\tilde{I}_{t+L-1}]$ as
\begin{equation*}
\begin{aligned}
E[\tilde{I}_{t+L-1}] &= \sum_{k=1}^{L}{\int_{\tilde{\xi}_{k-1}}^{\tilde{\xi}_k}(1-P\{x+\sum_{m=1}^{n}\tilde{X}_m \geq \tilde{\xi}_n,n=k,...,L\})dx}.   
\end{aligned}
\end{equation*}
Introducing the uniformly distributed random variables $W_k, k=1,...,L$,
\begin{equation*}
    W_k \overset{d}{=} U(0,\tilde{\xi}_k-\tilde{\xi}_{k-1}),
\end{equation*}
we find the following expression for $E[\tilde{I}_{t+L-1}]$,
\begin{equation*}
\begin{aligned}
E[\tilde{I}_{t+L-1}] &= \sum_{k=1}^{L}{(\tilde{\xi}_{k-1}}-{\tilde{\xi}_k})(1-P\{W_k+\sum_{m=1}^{n}\tilde{X}_m \geq \tilde{\xi}_n-\tilde{\xi}_{k-1},n=k,...,L\}).   
\end{aligned}
\end{equation*}
This equation can be rewritten with the help of the following (random) variables for $k=1,...,L$,
\begin{equation*}
\begin{aligned}
    \hat{X}_{k1} &= W_k+\sum_{m=1}^k\tilde{X}_m\\
    \hat{X}_{kn} &= \tilde{X}_{k+n-1}, n=2,...L+1-k\\
    \hat{\xi}_{kn} &= \tilde{\xi}_{k+n-1}-\tilde{\xi}_{k-1}, n=1,...,L+1-k.
\end{aligned}
\end{equation*}
Thus, we find an \textbf{exact} expression for $E[\tilde{I}_{t+L-1}]$,
\begin{equation}
\label{eqPIL}
\begin{aligned}
E[\tilde{I}_{t+L-1}] &= \sum_{k=1}^{L}{(\tilde{\xi}_{k-1}}-{\tilde{\xi}_k})(1-P\{\sum_{m=1}^{n}\hat{X}_{km} \geq \hat{\xi}_{kn},n=1,...,L+1-k\}).   
\end{aligned}
\end{equation}
This concludes this section in which we developed exact expressions for $P_{3,t+L}$ and $E[\tilde{I}_{t+L}]$.  

\section{Computational considerations}
\label{Compcons}
In this section we provide exact and approximate algorithms for computing the expressions for the $FP_3$-policy, the $PIL$-policy, and the constant order policy. These algorithms are used to determine the replenishment order $Q_t$ for a given state of the system at the start of period $t$. For other policies we implemented similar computational methods as proposed in \cite{zipkin2008old}, \cite{xin2021understanding}, and \cite{van2024projected}. These existing methods are discussed in section \ref{sec5}.

For the $FP_3$-policy we use bisection to solve for $Q_t$ in
\begin{equation}
\label{Q_P3}
\begin{aligned}
P\{\tilde{I}_{t+L}>0\} &= P_3.  
\end{aligned}
\end{equation}
For the $PIL$-policy it suffices to determine $E[\tilde{I}_{t+L-1}]$. Then we find $Q_t$ from
\begin{equation}
\label{Q_PIL}
\begin{aligned}
Q_t &= \tilde{I} - E[\tilde{I}_{t+L-1}],   
\end{aligned}
\end{equation}
with $\tilde{I}$ the target projected inventory level.

Finding the cost-optimal $P_3^*$ exploits the numerical finding that $\frac{E[T^2]}{E[T]}$ is increasing in $P_3$. Alternatively, we find the cost-optimal $P_3^*$ from assuming that the cost function is convex in $P_3$. The latter approach is relevant for discrete demand distributions, where Theorem \ref{lem:Optimality equation for lost-sales systems} does not hold. Finding the cost-optimal $\tilde{I}^*$ uses the convexity of the cost function under the $PIL$-policy (cf. \cite{van2024projected}). 

\subsection{Algorithms for non-stockout probability}
\label{algp3}
In subsection \ref{subsec41} we derived an exact backward recursion for the non-stockout probability. The recursion builds on the fact that $X_n$ is independent of $(Z_{n-1}-(\xi_{n-1}-\xi_{n-2}|Z_{n-1}>\xi_{n-1}-\xi_{n-2})$. This yields the following set of recursive equations,
\begin{equation*}
\begin{aligned}
E[Z_1] &= E[X_1],\\
E[Z_n] &= E[X_n]+E[(Z_{n-1}-(\xi_{n-1}-\xi_{n-2}|Z_{n-1}>\xi_{n-1}-\xi_{n-2})], n=2,...,L+1.\\
\sigma^2(Z_1) &= \sigma^2(X_1),\\
\sigma^2(Z_n) &= \sigma^2(X_n)+\sigma^2((Z_{n-1}-(\xi_{n-1}-\xi_{n-2}|Z_{n-1}>\xi_{n-1}-\xi_{n-2})), n=2,...,L+1.\\
\end{aligned}    
\end{equation*}
The first two moments of $(Z_{n-1}-(\xi_{n-1}-\xi_{n-2})|Z_{n-1}>\xi_{n-1}-\xi_{n-2})$ are easy to compute for probability distributions such as gamma, mixed Erlang, and shifted exponential (cf. \cite{tijms2003first}, Appendix \ref{Appendix_A}). Thus, we can find expressions for the first two moments of $Z_n, n=1,...,L+1$, and fit a distribution on these moments to obtain an approximation for $P\{Z_n>\xi_n-\xi_{n-1}\}, n=1,...,L+1$. It follows from equations \eqref{eqP3}, \eqref{eqYn}, and \eqref{eqZn} that   
\begin{equation}
\label{eqP32}
\begin{aligned}
P_{3,t+L}&=1-\prod_{n=1}^{L+1}P\{Z_n>\xi_n-\xi_{n-1}\}.
\end{aligned}    
\end{equation}

Equation \eqref{eqP32} gives an exact expression for the non-stockout probability $P_{3,t+L}$. Together with the approximations for $P\{Z_n>\xi_n-\xi_{n-1}\}, n=1,...,L+1$, from the two-moment recursion scheme, we find an approximation for $P_{3,t+L}$. In section \ref{sec5} we show that this approximation is remarkably accurate when assuming demand per period is mixed Erlang distributed and fitting mixed Erlang distributions on the first two moments of $Z_n, n=1,...,L+1.$ 
At every time $t$, the fixed non-stockout probability ($FP_3$) policy, solves $Q_t$ from equation \eqref{eqP32} with $P_{3,t+L}=P_3^*$. As $Q_t=\xi_1$ we choose some initial value of $Q_t$, subsequently apply the moment-iteration scheme to compute $P_{3,t+L}$ from equation \eqref{eqP32}, and compare $P_{3,t+L}$ and $P_3^*$. Thus, a bisection scheme yields $Q_t$.

We note here that the moment iteration scheme proposed in \cite{Dons1996} could be applied here, too. However, we found that this forward recursion, though more efficient in computing $Q_t$, yields less accurate results. I.e., when setting a target $P_3^*$, the simulated $P_3$ value of the backward recursion is closer to $P_3^*$ than the $P_3$ value of the forward recursion.
Furthermore, we stipulate that, when applying the $FP_3$ policy in practice, we need to calculate only a single $Q_t$ each period for each SKU. With current computational power, this is a matter of seconds for thousands of SKUs.

Where possible, we benchmark the quality of the approximations against exact expressions for $P_{3,t+L}$. For exact expressions for discrete demand distributions we can use the forward recursive scheme defined by equation \eqref{eqbal} and standard probability theory. In \cite{Dons1996} an exact expression is given for $P_{3,t+L}$ for pure Erlang distributions. However, the algorithm proposed is numerically intensive for large $L$ and low coefficients of variation of the demand. In \cite{van2024projected} an elegant and efficient algorithm is proposed for exact computation of $E[\tilde{I}_{t+L-1}]$ when the demand distribution is a mixed Erlang distribution with a single rate parameter. Instead of describing demand for items, they propose to consider the demand for exponential phases and describing the system state as the probability distribution of initial available phases to satisfy demand from, and a probability distribution of phases associated with each outstanding order. With a mixed Erlang distribution consisting of two pure Erlang distributions, demand for phases is represented by a two-point distribution. Thus, using standard laws of probability, we compute the distribution of the number of available phases at time $t+s$ from the probability distribution of the available phases at time $t+s-1$, the two-point demand distribution, and the probability distribution of phases added by replenishment order $Q_{t+s-L}$. The memorylessness of the exponential distribution enables to determine the  $E[\tilde{I}_{t+L-1}]$ from the inventory balance equation and the average amount lost each period. We can use the memorylessness property to argue that the probability of lost-sales in the original systems is equal to the probability of lost-sales in the "exponential phases system".
Fitting a Mixed Erlang distribution with a single rate parameter is straightforward. For sake of completeness and to correct errors in \cite{van2024projected} we present the expressions for the parameters of the Mixed Erlang distributions in appendix \ref{Appendix_A}.

The algorithm proposed in \cite{van2024projected} can also be applied to so-called shifted mixed Erlang distributions. Such distributions are the sum of a constant and a mixed Erlang distributed random variable. In our context the constant can be interpreted as a minimum demand per period. The special case of a shifted exponential distribution allows for two-moment fits as long as the standard deviation of demand is not exceeding the mean demand. By subtracting the constant from available inventory and each outstanding order, the system is equivalent to a lost-sales system with exponential demand. In the "exponential phases system" this implies that demand is constant equal to 1, which reduces the computational effort substantially. This is especially important for systems with low coefficient of variation, as in that case demand consists of many exponential phases, which blows up the relevant state space.
In appendix \ref{Appendix_A} we present the expressions for the parameters of the shifted exponential distribution.

\subsection{Algorithms for the PIL-policy}
To determine $Q_t$ for the PIL-policy we note that the probabilistic expressions on the right-hand side of equation \eqref{eqPIL} are identical to the expression for $P_{3,t+L}$ in equation \eqref{eqP3}. Thus, we can reuse our approximation scheme based on fitting mixtures of Erlang distributions to compute an \textbf{approximation} for $E[\tilde{I}_{t+L-1}]$. Using this scheme, it follows that computation of the approximation for $E[\tilde{I}_{t+L}]$ requires a factor of $\frac{L}{2}$ more CPU time.

As mentioned above, \cite{van2024projected} propose an exact algorithm for $E[\tilde{I}_{t+L-1}]$ for mixed Erlang distributions with a single rate parameter. Like for the expression for $P_{3,t+L}$, this exact algorithm can also be used to compute $E[\tilde{I}_{t+L-1}]$ for the case of shifted mixed Erlang distributions. For discrete demand distributions we can use again the recursive equation \eqref{eqbal} and standard probability theory.    

\subsection{Algorithms for the constant order policy}
\label{subsec_co}
The lost-sales system under a constant order policy is equivalent to a $G/D/1$ queue. \cite{Kok1989} gives an accurate approximation for the $G/D/1$ key performance characteristics, based on a moment-iteration method. In fact, this is the forward recursion applied in \cite{Dons1996}. We use the constant order policy as a benchmark for our Simulation Based Optimization approach discussed in section \ref{sec5}. As the constant order policy is not responsive to the system state, we expect that the time until stationarity is longest for this policy. Interestingly, under the assumption of shifted exponential demand, we can exploit the equivalence of the lost-sales system and the $G/D/1$ to find closed-form expressions for the optimal order quantity $Q^*$ and the minimum long-run average costs $E[C^*]$. We can transform the $SE/D/1$ queue to the $M/D/1$ queue for which we have closed-form expressions for the average waiting time and probability of waiting. This yields
\begin{equation}
\label{q_co}
\begin{aligned}
Q^* &=  E[D] \left(1-c_D\sqrt{\frac{h}{2p+h}}\right);
\end{aligned}    
\end{equation}
\begin{equation}
\label{ec_co}
\begin{aligned}
E[C^*] &=h \frac{\left(Q^*-(E[D]-\sigma(D)\right)^2}{2(E[D]-Q^*)} + p(E[D]-Q^*);
\end{aligned}    
\end{equation}
\begin{equation}
\label{p3_co}
\begin{aligned}
P_3^* &= 1- \sqrt{\frac{h}{2p+h}}.
\end{aligned}    
\end{equation}

Equations \eqref{q_co}, \eqref{ec_co}, and \eqref{p3_co} give insight how the constant replenishment order depends on the parameters $c_D$, $p$, and $h$. Quite remarkable is the fact that the optimal $P_3$ under the Constant Order policy does not depend on the coefficient of variation. In some of our numerical experiments in section \ref{sec5} we investigate the robustness of the optimal $P_3$ under various policies and different demand distributions. 

\section{Numerical experiments}
\label{sec5}
In this section we report our findings concerning the performance of the $FP_3$-policy in comparison to other policies proposed in the literature. Like in earlier papers on the lost-sales model, we use simulation-based optimization (SBO) to determine optimal policy parameters. We validate our SBO methodology by comparing the SBO results for the Constant Order (CO) policy with the exact closed-form expressions given by equations \eqref{q_co} and \eqref{ec_co} for the optimal order quantity and the long-run average costs, respectively. Next, we benchmark the $FP_3$ policy and the $PIL$ policy against the optimal policy, using the cases presented in \cite{zipkin2008old}. As the cases in \cite{zipkin2008old} concern specific cases only, we run an extensive discrete event simulation experiment, where we compare the performance of several policies over a wide range of parameter values, in particular including cases with long lead times and cases with very high penalty costs. This helps us to assess the performance of these policies under the asymptotic regimes for which the theoretical results derived in \cite{goldberg2016asymptotic} and \cite{huh2009asymptotic} hold. The policies compared are the Base Stock ($BS$-) policy, the $CO$-policy, the Capped Base Stock ($CBS$-) policy proposed in \cite{xin2021understanding}, the $PIL$-policy and the $FP_3$-policy proposed in this paper.

Note that the extensive simulation experiment also validates the SBO methodology: for the $FP_3$-policy and $PIL$-policy we determine optimal parameters that are both input and output in the simulation experiment. As the SBO methodology may be computationally intensive for exact algorithms, we investigate the performance of alternative algorithms applicable to the same first two moments of the demand distribution. This also allows us to investigate the sensitivity of the optimal policy parameters for different demand distributions with the same first two moments.    

\subsection{The simulation-based optimization (SBO) methodology}
\label{subsec_SBO}
Analysis of the lost-sales model under the policies proposed in literature does not yield analytical expressions for the long-run average costs, with the exception of the CO-policy (see section \ref{subsec_co}). Thus we need to resort to numerical approaches. We can find the optimal policy from policy iteration (cf. \cite{zipkin2008old}), but this is only possible for short lead times due to state space explosion as the lead time gets longer. Alternatively, authors use discrete event simulation to estimate the long-run average costs given the policy parameters and then use some optimization method, possibly exploiting properties such as convexity of the cost function, to determine the optimal policy parameters. We adopt the latter approach, as even the replenishment decision is implicit under the $FP_3$-policy, which is also the case for the $PIL$-policy.

In order to find the optimal policy parameters, we used a random number generator with random seeds to ensure that all simulations used the same realizations of the demand process. By doing so, we can prove convexity of the cost function under the $BS$-policy with respect to the base stock level $S$. Likewise, we can prove convexity of the cost function under the $CO$-policy with respect to the order quantity $Q$. As pointed out in \cite{xin2021understanding}, under the $CBS$-policy,  the cost function is not convex. However, we found that a nested optimization procedure where we assumed convexity of the cost function given a fixed maximum order quantity $Q_{max}$ with respect to the base stock level $S$, and convexity in the neighbourhood of the optimal maximum order quantity $Q_{max}^*$ yielded the optimal $CBS$-policy parameter values. We note here that we used a small grid size for $Q_{max}$ to determine the appropriate neighbourhood. We used the optimality equation \eqref{eq1} and bisection to determine the optimal $FP_3$-policy. Here we exploited the monotonicity of $\frac{E[T^2]}{E[T]}$ as a function of the penalty cost $p$ and thereby of $P_3$. We denote this policy by $FP_3-$opt in our results. As finding the optimal $P_3$ target value is computationally intensive, we also considered a heuristic, where the $P_3$ target value is determined by taking the $P_3$ value of the optimal $BS$-policy or optimal $CO$-policy, whichever performed best. We denote this policy by $FP_3-$heuristic in the results of our experiments in section \ref{subsec52}. As pointed out in \cite{van2024projected}, we used the convexity of the cost function with respect to the $PIL$ target value to determine the optimal target $PIL$. For the $FP_3$- and $PIL$-policy, we used both exact and approximative algorithms to determine each time unit the replenishment quantity.

After exploratory experimentation, we concluded that a simulation run length of 100,000 time units suffices to get accurate estimates of the long-run average costs under the various policies. As the $CO$-policy does not apply feedback using the system state, we expect that the time needed to reach stationarity is longest under this policy. We assumed a shifted exponential distribution (cf. appendix \ref{Appendix_A}) so that we can compare the SBO results with the exact results according to equations \eqref{q_co} and \eqref{ec_co}.

\begin{table}
    \centering
    \includegraphics[width=0.9\linewidth]{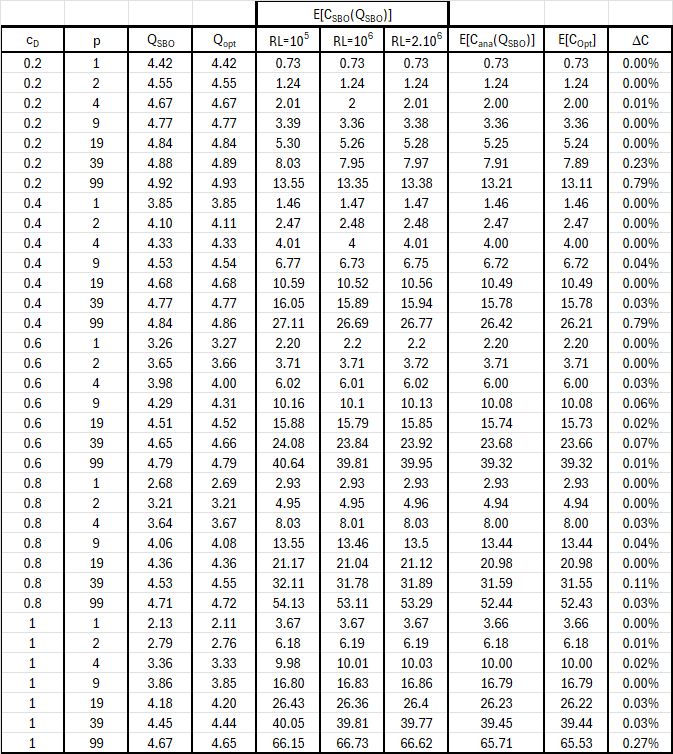}
    \caption{Results SBO methodology for the $CO$-policy}
    \label{Table_1}
\end{table}

As we used a run-length of $10^5$ time units for policy optimization, we tested the validity of the average costs performance for this run length by comparing against the average costs performance under the same policy for run lengths of $10^6$ and $2\cdot10^6$. The results of this experiment are presented in table \ref{Table_1}. We varied the coefficient of variation of demand $c_D$, and the penalty cost $p$. The holding cost $h$ is normalized at a value of 1. We can conclude that the policies determined by the SBO policy are close to optimal and that average costs do not change much when increasing the run length. We do see that under high penalty costs and thereby high order quantity, the $G/D/1$ utilization effect emerges, implying that even after $2\cdot10^6$ time units, the average costs for lost-sales system under the $CO$-policy has not converged to the long-run average costs.

\subsection{Benchmarking against the Zipkin cases}
\label{subsec_Zipkin}
In \cite{zipkin2008old}, for a set of cases, the long-run average costs under different policies are compared with the minimum long-run average costs under the optimal policy. The optimal policy is computed by a policy iteration algorithm with a discount factor of $0.995$. This set of cases has been the benchmark for most recent papers on the lost-sales model. Demand is Poisson and geometrically distributed with mean 5, $p$ is varied as 4, 9, 19, 39, $L$ is varied as 1, 2, 3, 4, $h$ equals 1. The parameter $c_D$ denotes the coefficient of variation of the demand per period. For both the Poisson distribution and the geometric distribution, $c_D$ is determined by the mean. In \cite{van2024projected} the $PIL$-policy is shown to be superior to other policies proposed in literature, so we compare the $FP_3$-policy and $PIL$-policy against the optimal policy. Left of the optimal costs column in table \ref{Table_2} we show the performance of the policies using the exact algorithms. We determined the policy parameter with a run-length of $10^5$, and evaluated it with a run-length of $10^5, 10^6$, and $2\cdot10^6$. We can conclude that the $FP_3$-policy outperforms the $PIL$-policy in costs and performs very close to optimal for all cases. To the right of the optimal costs column we show the performance of the forward and backward approximation algorithms for both policies. We conclude that the backward algorithm for the $FP_3$-policy performs almost identical to the exact algorithm, while the efficient forward algorithm for the $PIL$-policy performs almost identical to the exact algorithm. 

\begin{table}
    \centering
    \includegraphics[width=\linewidth]{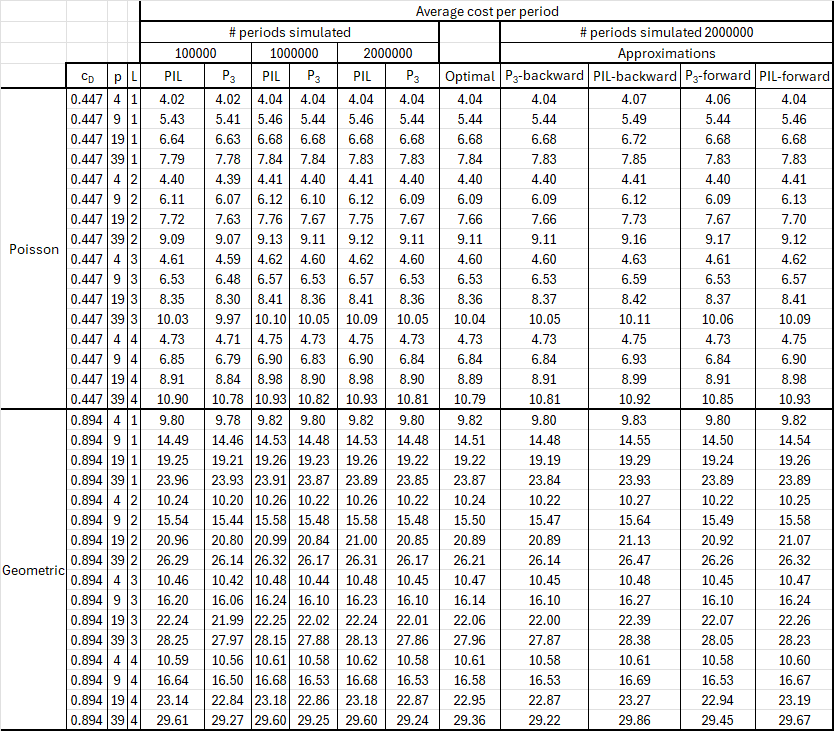}
    \caption{Comparison against the optimal policy}
    \label{Table_2}
\end{table}

We should make a remark about the observation that the $FP_3$-policy outperforms the optimal policy, even using a run-length of $2\cdot10^6$ time units. We noted above that the optimal policy in \cite{zipkin2008old} is computed with policy iteration with a discount factor of $0.995$. The difference in long-run average costs may stem from the use of the high discount factor instead of using a value iteration algorithm designed for long-run average costs minimization.

\subsection{Benchmarking against the Xin cases}
\label{subsec_Xin}
In \cite{xin2021understanding} the Zipkin cases are extended to $L=6, 8, 10$. In table \ref{Table_3} we compare the results of the $P_3$-policy and $PIL$-policy, derived from using an exact algorithm for discrete demand, with the results for the $CBS$-policy, $BS$-policy, and $CO$-policy. We may conclude that both the $PIL$- and $CBS$-policy perform much better than the $BS$- and $CO$-policy, but the $FP_3$-policy outperforms all policies especially for higher penalty costs. We also conclude that the SBO methodology reproduces the results in \cite{xin2021understanding}.\\

\begin{table}
    \centering
    \includegraphics[width=\linewidth]{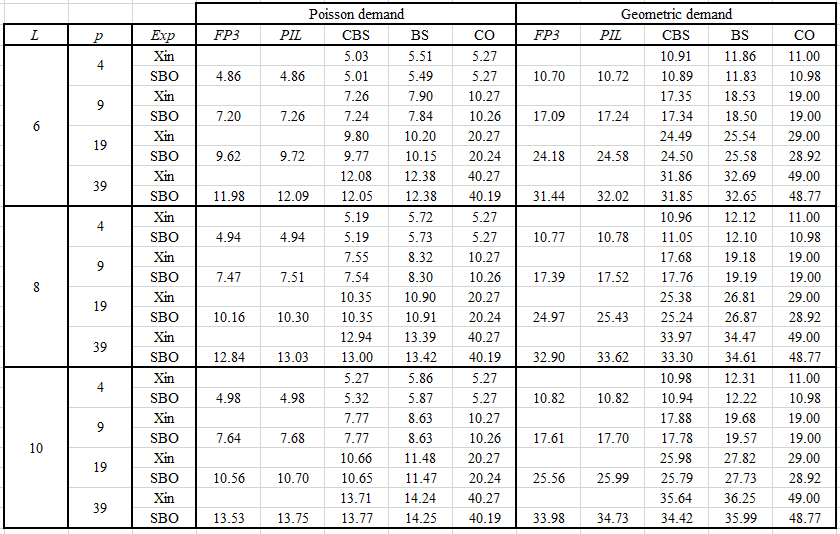}
    \caption{Comparison with Xin (2021) cases}
    \label{Table_3}
\end{table}

%\textcolor{red}{\textbf{De shifted exponential distribution is belangrijk instrument om efficient te rekenen. Blijkt dat de optimale $P_3$ voor ME en SE gelijk is. De kostenuitkomsten niet. Maar dat betekent dat een robuuste heuristiek voor $c_D$ is om met SE de optimale $P_3$ te bepalen en dan met ofwel de exacte ME of the backward recursion de kosten te bepalen. Voor $c_D>1$ is er niet zo'n efficiënte methode. Dan is de backward recursion eigenlijk de enige methode, die eventueel op 3 momenten gefit kan worden. Daarnaast is er de $P_3$-heuristiek, die ook goede resultaten geeft. Ik schrijf dit op, omdat hiermee de numerieke studie weer logisch wordt. }}

\subsection{Comparing policies in an extensive experiment}
\label{subsec52}
So far we focused on the performance of the $FP_3$-policy under discrete demand and for a limited set of model parameters. We developed approximations for the $FP_3$- and $PIL$-policy to apply them to a wide range of practical situations under the assumption of continuous demand per period. Therefore, we ran an extensive discrete event simulation experiment. In table \ref{tab: parameters} we present the model parameter combinations used. We considered very high values of $p$ and $L$ to explore the asymptotic optimality of the base stock policy and constant order policy, respectively, in comparison with the $FP_3$ and $PIL$ policies. We assumed demand to be a mixture of Erlang-(k-1) and Erlang-k distributions when $c_D \leq 1$, and hyperexponentially distributed with gamma moments when $c_D > 1$ (cf. appendix \ref{Appendix_A}). We note here that a discrete event simulation experiment on a model with as little parameters as the periodic review lost-sales model can be conclusive on the performance of policies, even without formal mathematical analysis.

\begin{table}[ht]
    \centering
    \caption{Design of discrete event simulation experiment}
    \begin{tabular}{ll}\\
    
    \hline
        $p$ & 4, 9, 19, 49, 99, 199\\
        $h$ & 1\\
        $L$ & 1, 2, 4, 8, 16, 32, 64\\
        $c_D$ & 0.25, 0.5, 0.75, 1, 1.25, 1.5, 1.75, 2\\
        \hline
    \end{tabular}
    \label{tab: parameters}
\end{table}

\begin{figure}[h!]
  \centering
     \includegraphics[scale=0.9]{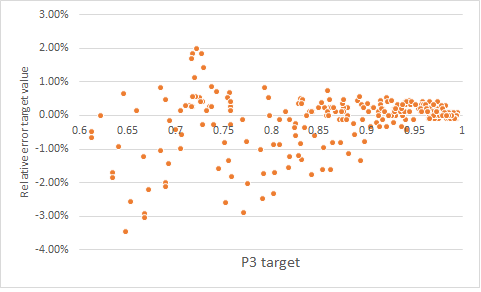}
  \caption{Accuracy of $P_3$ approximation}
  \label{fig:1}
\end{figure}

The main findings from our experiment are as follows:

1. The backward approximation of $P_3$ is very accurate for all parameters (cf. figure \ref{fig:1}). Errors are highest for low value of the penalty cost $p$, implying a low target value for $P_3$.\\

\begin{figure}[h!]
  \centering
     \includegraphics[scale=0.9]{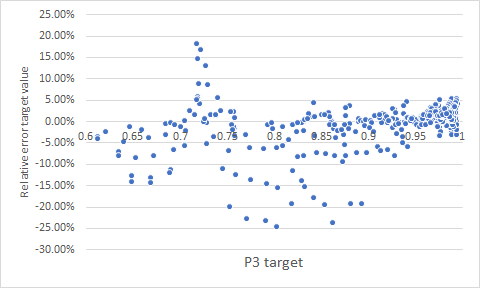}
  \caption{Accuracy of $PIL$ approximation as a function of the $FP_3$-target value}
  \label{fig:2}
\end{figure}

\begin{figure}[h!]
  \centering
     \includegraphics[scale=0.9]{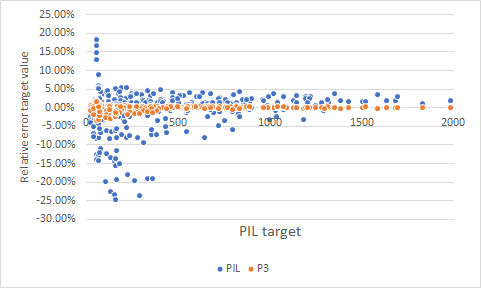}
  \caption{Accuracy of $PIL$ approximation as a function of the $PIL$-target value}
  \label{fig:3}
\end{figure}

2. The approximation of $E[I_{t+L}]$ is less robust (cf. figures \ref{fig:2} and \ref{fig:3}). In about $50\%$ of the cases the resulting average inventory from the simulation deviates more than $2\%$  from the target $PIL$. This can be explained by the higher relative errors in the $P_3$ approximation for low values of $P_3$ and equation \eqref{eqPIL}. In this equation we multiply the approximation for $P_3$ with an order quantity and sum $L$ such terms. If the error in these terms are all in the same direction, these errors add up. Errors are largest for cases where the $P_3$ target level is below 85\%. \\

\begin{figure}[h!]
  \centering
     \includegraphics[scale=0.9]{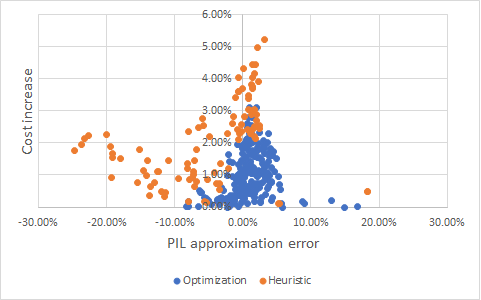}
  \caption{Impact of $PIL$ approximation error on performance approximate $PIL$-policy}
  \label{fig:12}
\end{figure}

\begin{figure}[h!]
  \centering
     \includegraphics[scale=0.9]{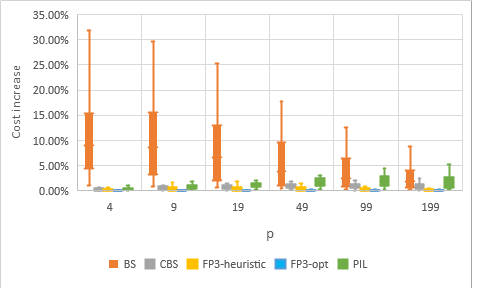}
  \caption{Cost increase with respect to best policy as a function of $p$, $CO$-policy excluded}
  \label{fig:8}
\end{figure}

\begin{figure}[h!]
  \centering
     \includegraphics[scale=0.9]{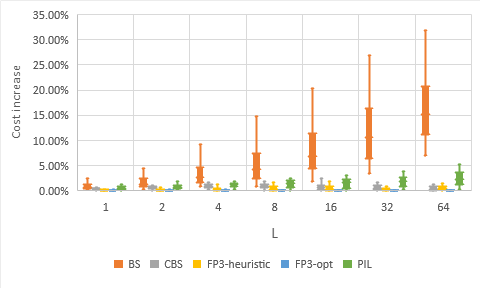}
  \caption{Cost increase with respect to best policy as a function of $L$, $CO$-policy excluded}
  \label{fig:9}
\end{figure}

3. We took a closer look at the impact of the (in)accuracy of the $PIL$-target approximation on the cost increase under the approximate $PIL$-policy. Figure \ref{fig:12} shows that when the approximation error is (very) high, the increase in cost is only limited, below 2\%. When the approximation error is small we see costs increases can be 2\%-5\%. From figures \ref{fig:8} and \ref{fig:9} we see that the $PIL$-policy performs worse for longer lead times and higher penalty costs. We see that $PIL$-target approximation errors are higher for the cases where we applied our heuristic, i.e. cases with $L=32, 64$, using the average stock at the start of a period from the $FP_3-opt$-policy. Still, we see that even then cost increases are limited for cases with high approximation errors, whereas costs increases can be high for cases with accurate approximations, typically when penalty costs are high. However, we should consider using the exact algorithm proposed in \cite{van2024projected} for the $PIL$-policy, which we extended to shifted ME distribution, when this is computationally feasible.    

\begin{figure}[h!]
  \centering
     \includegraphics[scale=0.9]{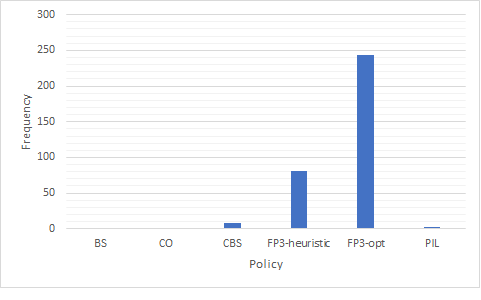}
  \caption{Best policy}
  \label{fig:4}
\end{figure}

4. Despite the approximation errors, the $FP_3$-policy is best in $97\%$ of all cases (cf. figure \ref{fig:4}) using the SBO methodology. The $BS$-policy was never best. The $CO$-policy was best in only 1 of the 336 cases. The $CBS$-policy was optimal in $2\%$ (8) of the cases. The $PIL$-policy was best in 2 cases, though, as stated above, this may be caused by the approximation error. The $FP_3$-heuristic performed best in $24\%$ of the cases, which we expect to be a result of the approximation error in $P_{3,t+L}$.
\newline

\begin{figure}[h!]
  \centering
     \includegraphics[scale=0.9]{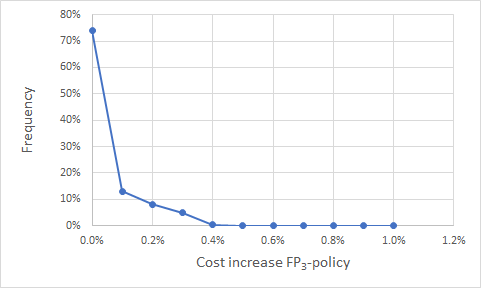}
  \caption{$FP_3$ policy against best policy}
  \label{fig:5}
\end{figure}

5. The $FP_3$-policy was never more than 0.4\% more expensive than the best policy (cf. figure \ref{fig:5}). When comparing with the analytically computed costs under the $CO$-policy, the maximum cost increase was $1.3\%$.
\newline

\begin{figure}[h!]
  \centering
     \includegraphics[scale=0.9]{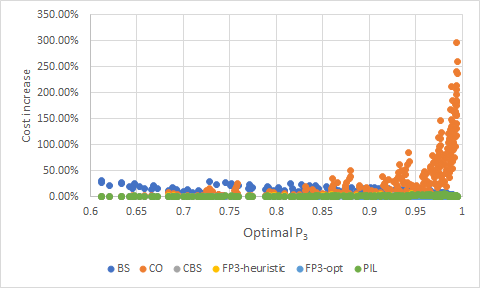}
  \caption{Cost increase with respect to best policy}
  \label{fig:6}
\end{figure}

6. The cost increase of the $CO$-policy can be dramatic compared to the best policy (cf. figure \ref{fig:6}). In particular for higher penalty costs, i.e for higher values of the target $P_3$ for the optimal $FP_3$-policy, cost increase by more than $50\%$.\\

\begin{figure}[h!]
  \centering
     \includegraphics[scale=0.9]{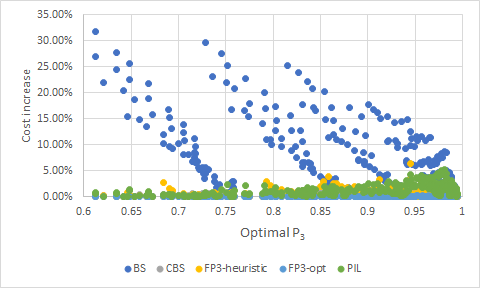}
  \caption{Cost increase with respect to best policy, $CO$-policy excluded}
  \label{fig:7}
\end{figure}

7. If we exclude the $CO$-policy from our comparison, we conclude that the BS policy is performing much worse than the other four policies, i.e. the $CBS$-policy, the two $FP_3$-policies and the $PIL$-policy (cf. figure \ref{fig:7}).\\

8. Exploring the performance of the policies further we find, as expected, that the $BS$-policy performs well for high penalty costs (cf. figure \ref{fig:8}), but not better than the $CBS$-policy and the $FP_3$-policy. However, when lead times are high, the asymptotic optimality of the $BS$-policy for high penalty costs cannot be seen (cf. figure \ref{fig:9}). The $CBS$-policy and $FP_3$-policy are both very robust, performing very well in comparison with the best policy for all cases.\\

\begin{figure}[h!]
  \centering
     \includegraphics[scale=0.9]{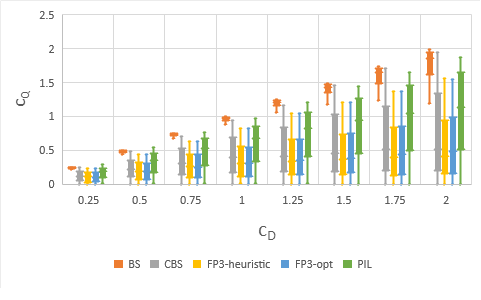}
  \caption{Coefficient of variation of replenishment process as a function of $c_D$}
  \label{fig:10}
\end{figure}

\begin{figure}[h!]
  \centering
     \includegraphics[scale=0.9]{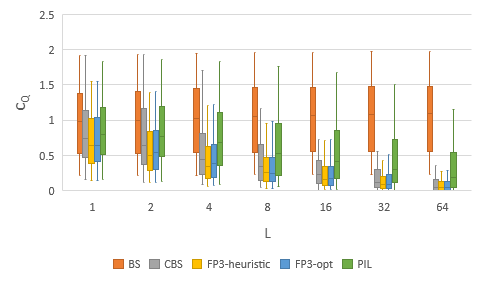}
  \caption{Coefficient of variation of replenishment process as a function of $L$}
  \label{fig:11}
\end{figure}

9. The coefficient of variation of the replenishment order under the $FP_3$-policy and $CBS$-policy is substantially lower than the coefficient of variation of the demand per period (cf. figure \ref{fig:10}). The replenishment order volatility approaches 0 as lead times get longer, thereby showing that the $FP_3$-policy and $CBS$-policy converge to the $CO$-policy as lead times increase. We see a similar effect under the $PIL$-policy, albeit much less pronounced (cf. figure \ref{fig:11}). From a practical point of view, this may be the most important conclusion from our numerical study.\\

The extensive experiment shows that the promising results regarding the performance of the $FP_3$-policy presented in subsections \ref{subsec_Zipkin} and \ref{subsec_Xin} extend to a wide range of model parameters. Thus, it is worthwhile to explore some more aspects of the lost-sales model under the $FP_3$-policy.

\subsection{Sensitivity of optimal policy parameters}
\label{subsec_sens}
In the numerical experiment we assumed that demand is distributed according to a specific mixture or Erlang distributions or hyperexponential distributions. Although these distributions have proven to be highly valuable when modelling real-life demand processes (cf. \cite{Kok2018}), we investigate the impact of this demand process assumption. Towards this end we compare the performance under this demand process assumption with shifted exponential demand ($c_D \leq 1$) and demand distributed according to a mixture of exponential and Erlang-k demand ($c_D>1$).  We varied $c_D$ as 0.25, 0.5, 0.75, 1, 1.25, 1.5, 1.75, and 2. We varied $p$ as 4, 19, and 99. We varied $L$ as 1, 2, 4, 6, 8, and 10. For $c_D \leq 1$, ME refers to the mixture of Erlang-(k-1) and Erlang-k, while for $c_D>1$, ME refers to the mixture of exponential and Erlang-k. SE refers to shifted exponential and HY refers to hyperexponential with gamma moments. For the SE and ME cases we used exact algorithms, for the HY cases we used the backward recursion algorithm. The results of the demand distribution sensitivity study are shown in Table \ref{Table_PDF}. 

\begin{table}
    \centering
    \includegraphics[width=\linewidth]{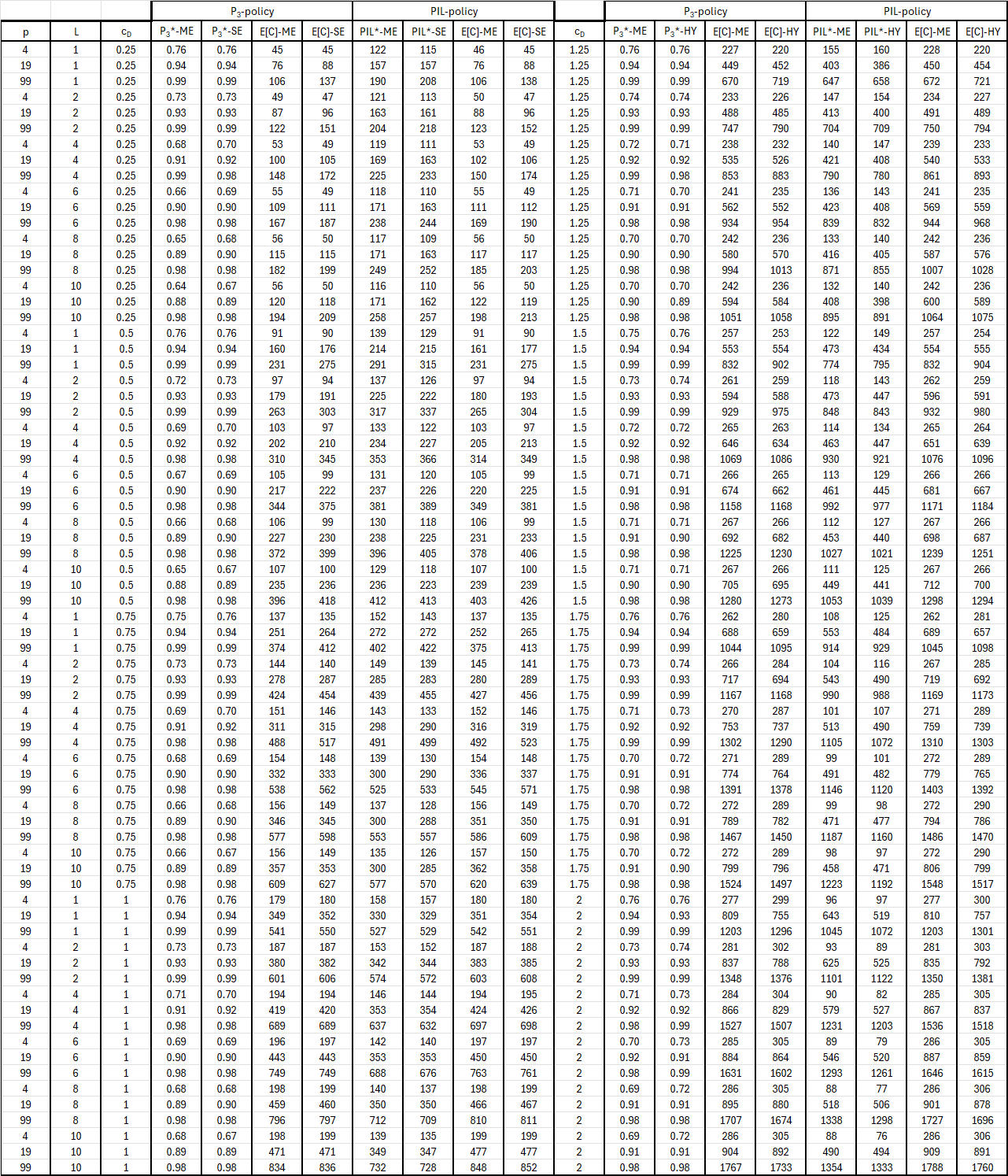}
    \caption{Policy parameter and cost sensitivity to the demand distribution}
    \label{Table_PDF}
\end{table}

The optimal $P_3$ target level is insensitive to the distribution over all parameter values. The average costs are more sensitive to the distribution assumption as lead times get longer and penalty costs get higher. As expected, the average costs are more sensitive to the distribution as $c_D$ gets closer to 0 and as it gets closer to 2. The optimal $PIL$ target level is sensitive to the distribution as $p$ increases: the difference between the optimal $PIL$ target values may change in sign. Note that the $FP_3$-policy outperforms the $PIL$-policy for all cases with $c_D \leq 1$ and for 71 out of 72 cases for the ME distribution and 67 out of 72 cases for the HY distribution. The $PIL$-policy only outperformed the $FP_3$-policy for $c_D=2$ with a very small percentage. For the specific case of $c_D=2$, $p=19$, and $L=2$ and the ME distribution, we applied the SBO-methodology with longer run-lengths for optimization to find that the $PIL$-policy outperforms the $FP_3$-policy. This provides clear evidence that, while the $FP_3$-policy overall performs better than the other policies, it is not the optimal policy.

\subsection{Practical application of the SBO methodology}
\label{subsec_practice}
In view of the above, we propose a procedure to apply the $FP_3$-policy in practice. While we so far used an optimization run-length of $10^5$, we found that an optimization run-length of $10^4$ suffices, which reduces computation times by a factor of 10. For the performance evaluation of the policy found, we recommend a run-length of $10^6$. Given the insensitivity of the target $P_3$, we recommend to fit the shifted exponential distribution (cf. subsection \ref{dis_SE}) to the first two moments of the demand distribution for $c_D \leq 1$. For $c_D>1$ we recommend to use the hyperexponential distribution with the first three moments equal to those of the ganmma distribution (cf. subsection \ref{dis_hyp}). For convenience we have applied this procedure to generate the optimal $P_3$ target value for a wide range of parameters. The results can be found in the look-up table \ref{Table_A1} in appendix \ref{Appendix_B}. From this table we find that for low values of the penalty costs $p$, the target $P_3$-value is constant. Further investigation showed that for low values of the penalty costs, the $FP_3$-policy behaves like a Constant Order policy. This is confirmed by the long-run average costs under the $FP_3$-policy as presented in table \ref{Table_A2}.

\section{Conclusions}
\label{sec6}
In this paper we developed new results for the classical periodic review lost-sales model with linear holding and penalty costs. We investigated the performance of the $FP_3$-policy introduced by \cite{Dons1996}, deriving an exact backward recursion algorithm to compute the replenishment quantities under this policy. The $FP_3$-policy is optimal for the periodic review backlog model with linear holding and penalty costs. 
For the lost-sales model, we formulated a new optimality equation that should hold under the optimal policy, despite the fact that this policy is to-date unknown. This optimality equation applies to other policies as well, such as the $CO$-policy and the $FP_3$-policy. In the latter case, the optimality equation enables us to find the optimal $FP_3$-policy using discrete event simulation. 

We calibrated a simulation-based optimization (SBO) methodology based on comparing the SBO results for the Constant Order policy with closed-form expressions for the optimal order quantity and long-run average costs for shifted exponential demand. We extended the exact analysis in \cite{van2024projected} of the $PIL$-policy for Mixed Erlang demand to the exact analysis of the $FP_3$-policy and $PIL$-policy for shifted Mixed Erlang demand. We showed that the optimal $FP_3$-policy is very close to optimal for the discrete demand cases in \cite{zipkin2008old}. The results from our extensive numerical experiment show the robustness of the $FP_3$-policy proposed: the efficient approximation scheme results in a policy that outperforms policies for which asymptotic optimality has been proven formally. The results correspond with our proof that the $PF_3$-policy behaves like those asymptotically optimal policies for high penalty costs and long lead times.  

The various approximation algorithms to compute the replenishment quantity under the $FP_3$-policy all yield an excellent cost performance. However, the target $P_3$ value may deviate from the actual $P_3$ value. The target $P_3$ value from the backward recursion algorithm is very close to the actual $P_3$-value, and is therefore the preferred method when exact algorithms are inefficient.
In line with the observation in \cite{Dons1996}, we find a reverse Bullwhip effect for the lost-sales model: a good replenishment policy reduces the volatility in product demand towards the supplier of the product. The positive impact on the upstream supply chain of the volatility reduction found (cf. figures \ref{fig:10} and \ref{fig:11}) can not be overstated. Lower replenishment order volatility reduces the need for upstream inventory and reduces production costs.
Both the $CO$-policy and $BS$-policy are not recommended in a real-life setting, as cost increases compared to the costs of the $FP_3$-policy can be substantial, in particular for the $CO$-policy.
Our experiments provide further evidence of the excellent cost performance of the $CBS$-policy proposed by \cite{xin2021understanding}. It is essential to cap replenishment quantities in lost-sales systems with lead times of multiple time units. Apparently, this is a built-in feature of the $FP_3$-policy. In real-life situations, one can assume a target $P_3$ based on tacit knowledge from sales planners, whereafter the policy starts working immediately and caps when needed.
Our SBO methodology can be applied in practice. In case $c_D \leq 1$, fitting a Shifted Exponential distribution to the first two moments of the demand distribution and using the exact algorithm yields an efficient procedure with robust performance. In case $c_D > 1$, fitting a hyperexponential distribution to the first two or three moments and using the backward recursion algorithm, yields a procedure with acceptable efficiency and robust performance.

Obviously, there is a need for further research. The lost-sales model is a building block of other models, such as transshipments and multiple suppliers, where demand for an item is fulfilled by another source than originally planned for. It may be that the $FP_3$-policy can be a building block of policies for such systems. 
The approximation scheme does not assume demand stationarity, only demand independence in different periods. In a retail environment, demand is non-stationary, as demand on Mondays is typically lower than demand on Fridays and Saturdays. Assuming continuous demand, as we have done in this paper, we can still ensure that each period the target value can be met with a non-negative replenishment. Simulation experiments based on empirical retail data can show the impact of the application of the $FP_3$-policy.
As a closing remark, we again advocate the use of the non-stockout probability as the key inventory management performance indicator. It is easy to compute from standard available data, and setting the right target value yields a policy for both the backlog and lost-sales model that is likely to be (close-to-)optimal. \\

\newpage

\begin{APPENDICES}
\section{Two-moment distribution fits}
\label{Appendix_A}
For sake of self-containedness we provide the expressions for the demand distributions used to compute the replenishment quantities and to generate the demand in the discrete event simulation. The expressions below correct some typos in \cite{van2024projected}.

\subsection{Shifted exponential distribution}
\bigskip
\begin{equation*}
\label{dis_SE}
\begin{aligned}
d_0 &= (1-\sqrt{c_D})\\
\mu &= \frac{1}{\sqrt{c_D}E[D]}\\
F_D (x) &= 
    \begin{cases}
    0 & \text{if $x \leq d_0$}\\
    1-e^{-\mu(x-d_0)} & \text{if $x > d_0$}\\
    \end{cases}
\end{aligned}    
\end{equation*}
\\

\subsection{Mixture of Erlang-k and Erlang-(k-1)}
\bigskip
\begin{equation*}
\label{dis_Ek-1k}
\begin{aligned}
k &= \left\lfloor{\frac{1}{c_D^2}}\right\rfloor+1\\
q &= \frac{(kc^2-\sqrt{k(1+c_D^2)-k^2c_D^2})}{1+c_D^2}\\
\mu &= \frac{k-p}{E[D]}\\
F_D(x) &=1-q\sum_{i=0}^{k-2}{e^{-\mu x}\frac{x^i}{i!}}-(1-q)\sum_{i=0}^{k-1}{e^{-\mu x}\frac{x^i}{i!}}
\end{aligned}    
\end{equation*}
\\
\subsection{Mixture of Erlang-1 and Erlang-k}
\bigskip
\begin{equation*}
\label{dis_E1k}
\begin{aligned}
k &= \left\lfloor{2c_D^2+\sqrt{4(c_D^2-1)}}\right\rfloor+1\\
q &= \frac{2kc_D^2+k-2-\sqrt{k^2+4-4kc_D^2}}{2(k-1)(1+c_D^2)}\\
\mu &= \frac{q+k(1-p)}{E[D]}\\
F_D(x) &= 1-qe^{-\mu x}-(1-q)\sum_{i=0}^{k-1}{e^{-\mu x}\frac{x^i}{i!}}\\
\end{aligned}    
\end{equation*}
\\
\subsection{Hyperexponential distribution with gamma moments}
\bigskip
\begin{equation*}
\label{dis_hyp}
\begin{aligned}
\mu_1 &= \frac{2}{E[D]}\left(1+\sqrt{\frac{c_D^2-0.5}{c_D^2+1}}\right)\\
\mu_2 &= \frac{4}{E[D]}- \mu_1\\
q &= \mu_1\frac{\mu_2E[D]-1}{\mu_2-\mu_1}\\
F_D(x) &= 1 - qe^{-\mu_1 x}-(1-q)e^{-\mu_2 x}\\
\end{aligned}    
\end{equation*}
\\

\newpage
\section{Optimal parameters $FP_3-policy$}
\label{Appendix_B}
In this appendix we provide a look-up table for the optimal $P_3$ target value under a $FP_3$-policy. We varied the coefficient of variation of the demand $c_D$, the penalty costs $p$ and the lead time $L$. For $c_D \leq 1$ we assumed shifted exponential demand and for $c_D>1$ we assumed hyperexponential demand with the first three moments equal to those of the gamma distribution. For values of $c_D$, $p$, and $L$ not in the table, we propose to use some interpolation procedure.\\

\begin{table}
%    \centering
    \includegraphics[height=\linewidth, angle =90]{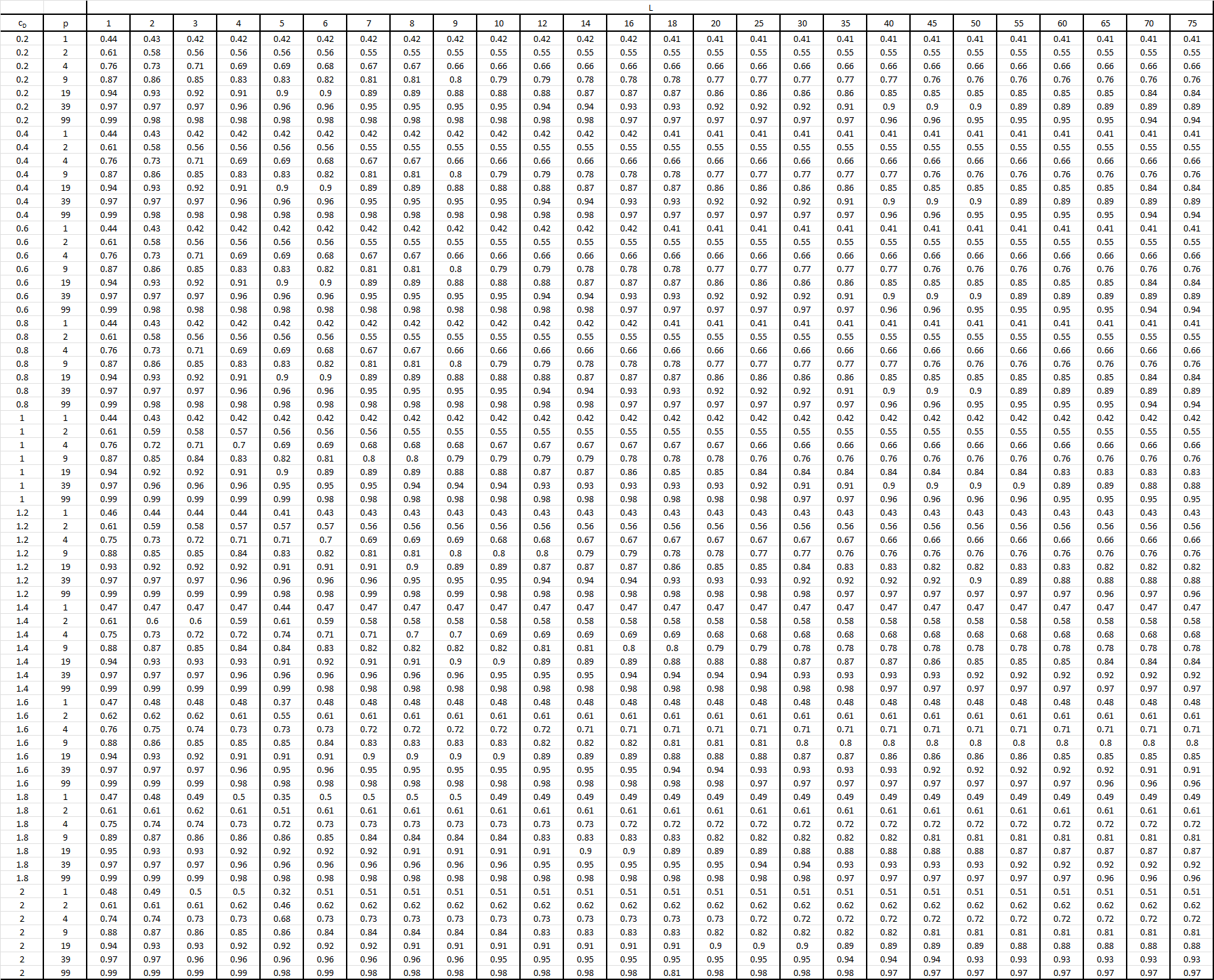}
    \caption{Optimal $P_3$ target value}
    \label{Table_A1}
\end{table}

\newpage

The long-run average costs under the optimal $P_3$ target value are presented in Table \ref{Table_A2}.

\begin{table}
%    \centering
    \includegraphics[height=\linewidth, angle =90]{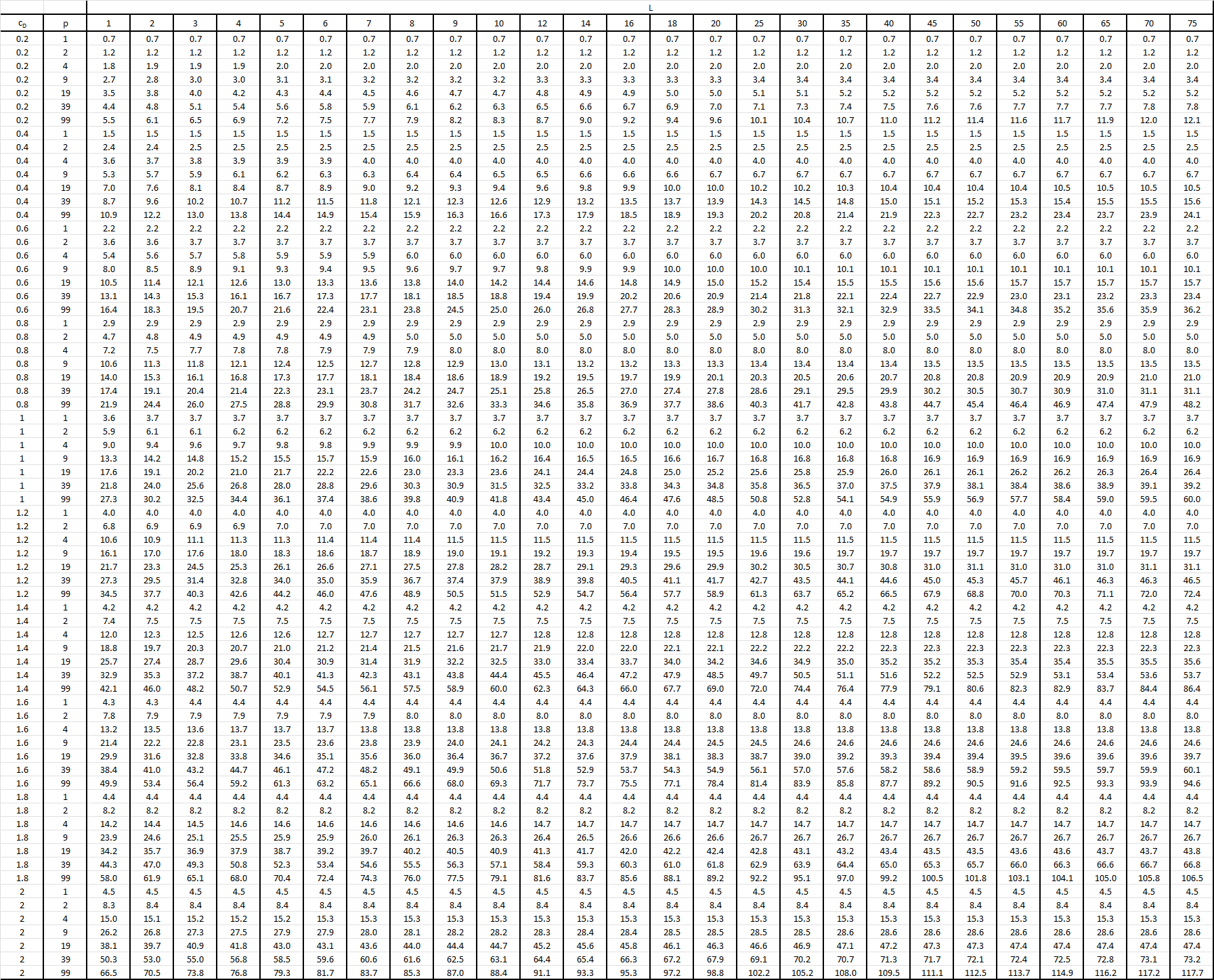}
    \caption{Long-run average cost under optimal $P_3$ target value}
    \label{Table_A2}
\end{table}

\end{APPENDICES}

% \begin{APPENDIX}{Comparing simulation and analysis of the CO policy}
% \label{Appendix_A}

% We noted that the lost-sales model under the $CO$-policy is equivalent to a G/D/1 queue. It is well-known that queueing systems under heavy load, higher than 90\%, say, need millions of arriving customers before the system reaches its stationary state. In our SBO methodology we decided to run 100,000 time units to determine the long-run average costs for a replenishment policy. In the table below we show that indeed when penalty costs increase, implying a higher fill rate, a higher order quantity, and thereby a higher utilization of the equivalent G/D/1 queue, the difference in SBO long-run costs and analytically derived long-run costs using the moment-iteration method from \cite{Kok1989} increases. Interestingly, the values of the optimal order quantities do not differ much. Note that, like in our extensive experiment, we assume Mixed Erlang distributed demand.\\

% ========  Insert Figure 13 about here ========\\

% \end{APPENDIX}

% \newpage

% \begin{APPENDIX}{Two-moment distribution fits}
% \label{Appendix_B}
    
% \end{APPENDIX}

% \newpage

% \begin{APPENDIX}{Optimal parameters $FP_3-policy$}
% \label{Appendix_C}
    
% \end{APPENDIX}
\newpage
\ACKNOWLEDGMENT{Herewith we would like to thank Karel van Donselaar for initiating our work on the lost-sales model after Ton's return to academia in 1992. Ever since, from time to time, Karel and Ton teamed up to work on developing algorithms for real-world problems, and get them implemented in practice. At the time they concluded that base stock policies are no good for lost-sales models and after this paper, our conclusion still stands. We would like to thank Joachim Arts and Willem van Jaarsveld for the many stimulating discussions on policies for the lost-sales model and almost forcing us to complete this paper.}
\newpage
\bibliographystyle{informs2014} % outcomment this and next line in Case 1
\bibliography{cas-refs} % if more than one, comma separated

\newpage
% Appendix here
% Options are (1) APPENDIX (with or without general title) or
%             (2) APPENDICES (if it has more than one unrelated sections)
% Outcomment the appropriate case if necessary
%
% \begin{APPENDIX}{<Title of the Appendix>}
% \end{APPENDIX}
%
%   or
%
% \begin{APPENDICES}
% \section{<Title of Section A>}
% \section{<Title of Section B>}
% etc
% \end{APPENDICES}

%%

% Acknowledgments here

% References here (outcomment the appropriate case)

% CASE 1: BiBTeX used to constantly update the references
%   (while the paper is being written).

% CASE 2: BiBTeX used to generate mypaper.bbl (to be further fine tuned)
%\input{mypaper.bbl} % outcomment this line in Case 2

%If you don't use BiBTex, you can manually itemize references as shown below.

%%%%%%%%%%%%%%%%%
\end{document}